\documentclass[leqno,12pt]{article}
\usepackage{amsmath,amssymb,color,a4wide}
\usepackage{MnSymbol}
\usepackage{tikz-cd}
\usepackage{xr-hyper}
\usepackage{hyperref}

\usepackage[curve,matrix,arrow,cmtip]{xy}
\usepackage{verbatim}
\NoComputerModernTips

\newdir^{ (}{{}*!/-3pt/\dir^{(}}
\newdir_{ (}{{}*!/-3pt/\dir_{(}}

\def\charact{\mathop{\rm char}\nolimits}
\def\BAfF{{\BA^{\rm f}_F}}
\def\bigtimes{\mathop{\raisebox{-2pt}{\huge $\times$\kern-1pt}}}

\def\BGacomma#1{{\BG_{\mathrm{a},#1}}}
\def\hatA{\smash{\hat{A}}}

\def\MM{{\mathcal M}}

\def\OM{\mathchoice
  {\rlap{\kern3.2pt$\overline{\phantom{L}}$}M}
  {\rlap{\kern3.2pt$\overline{\phantom{L}}$}M}
  {\rlap{\kern2.4pt$\scriptstyle\overline{\phantom{L}}$}M}
  {\rlap{\kern1.8pt$\scriptscriptstyle\overline{\phantom{L}}$}M}}

\let\le\leqslant
\let\ge\geqslant
\let\leq\leqslant
\let\geq\geqslant

 \def\Aut{\mathop{\rm Aut}\nolimits}

 \def\Cl{\mathop{\rm Cl}\nolimits}
 
 \def\Proj{\mathop{\rm Proj}\nolimits}
 \def\Spec{\mathop{\rm Spec}\nolimits}
 
 \def\deg{\mathop{\rm deg}\nolimits}
 
 \def\mod{\mathop{\rm mod}\nolimits}

\def\diag{\mathop{\rm diag}\nolimits}
\def\ord{\mathop{\rm ord}\nolimits}

\def\GL{\mathop{\rm GL}\nolimits}
\def\SL{\mathop{\rm SL}\nolimits}

\def\Id{{\rm Id}}

\let\phi\varphi
\let\theta\vartheta
\let\epsilon\varepsilon
\let\setminus\smallsetminus

\newtheorem{Thm}{Theorem}[section]
\newtheorem{Prop}[Thm]{Proposition}
\newtheorem{Lem}[Thm]{Lemma}

\newtheorem{Cor}[Thm]{Corollary}

\newtheorem{Def}[Thm]{Definition}

\newtheorem{Rem}[Thm]{Remark}
\newtheorem{Ex}[Thm]{Example}

\newtheorem{Ass}[Thm]{Assumption}

\def\UseTheoremCounterForNextEquation{\setcounter{equation}{\value{Thm}}\addtocounter{Thm}{1}}

\def\qed{{\hskip0pt\unskip\unskip\nobreak\hfil\penalty50
          \hskip1em\hbox{}\nobreak\hfil
           {$\square$}
          \parfillskip=0pt\finalhyphendemerits=0
          \par}\medskip}

\newenvironment{Proof}
               {\noindent{\bf Proof.}\ }
               {\qed}
							
               {\noindent{\bf Proof Sketch.}\ }
               {\qed}

               {\noindent{\bf Proof of #1.}\ }
               {\qed}

\newcommand{\BA}{{\mathbb{A}}}

\newcommand{\BC}{{\mathbb{C}}}

\newcommand{\BF}{{\mathbb{F}}}
\newcommand{\BG}{{\mathbb{G}}}

\newcommand{\BZ}{{\mathbb{Z}}}

\newcommand{\Fa}{{\mathfrak{a}}}

\newcommand{\Fp}{{\mathfrak{p}}}

\newcommand{\CD}{{\cal D}}

\newcommand{\CH}{{\cal H}}

\newcommand{\CL}{{\cal L}}
\newcommand{\CM}{{\cal M}}

\newcommand{\CO}{{\cal O}}

\newcommand{\CR}{{\cal R}}
\newcommand{\CS}{{\cal S}}

\newcommand{\CW}{{\cal W}}

\newbox\mybox
\def\arrover#1{\mathrel{
       \setbox\mybox=\hbox spread 1.4em
              {\hfil$\scriptstyle#1$\hfil}
       \vbox{\offinterlineskip\copy\mybox
             \hbox to\wd\mybox{\rightarrowfill}}}}

\def\larrover#1{\mathrel{
       \setbox\mybox=\hbox spread 1.4em
              {\hfil$\scriptstyle#1\vphantom{g}$\hfil}
       \vbox{\offinterlineskip\copy\mybox
             \hbox to\wd\mybox{\leftarrowfill}}}}

\def\ontoover#1{\mathrel{
       \setbox\mybox=\hbox spread 1.4em
              {\hfil$\scriptstyle#1\vphantom{g}$\hfil}
       \vbox{\offinterlineskip\copy\mybox
             \hbox to\wd\mybox{\rightarrowfill\hskip-2.8mm
                               $\rightarrow$}}}}
\def\leftontoover#1{\mathrel{
       \setbox\mybox=\hbox spread 1.4em
              {\hfil$\scriptstyle#1\vphantom{g}$\hfil}
       \vbox{\offinterlineskip\copy\mybox
             \hbox to\wd\mybox{$\leftarrow$\hskip-2.8mm
                               \leftarrowfill}}}}
\let\longto\longrightarrow

\def\Cinf{{\BC}_\infty}
\def\Finf{F_\infty}

\def\Bigskip{\bigskip\bigskip}

\newbox\dotDdbox
\newbox\dotDtbox
\newbox\dotDsbox
\newbox\dotDssbox
\setbox\dotDdbox=\vbox{\hbox{\kern3.5pt\bf.}\vskip-4.5pt\hbox{$D$}}
\setbox\dotDtbox=\vbox{\hbox{\kern3.5pt\bf.}\vskip-4.5pt\hbox{$D$}}
\setbox\dotDsbox=\vbox{\hbox{\kern2.1pt.}
                       \vskip-6.6pt\hbox{$\scriptstyle D$}}
\setbox\dotDssbox=\vbox{\hbox{\kern1.6pt.}
                        \vskip-8.1pt\hbox{$\scriptscriptstyle D$}}

\newcommand{\Theorem}{Theorem }

\newcommand{\PropositionCite}{Prop.\,}

\newcommand{\TheoremCite}{Thm.\,}
\newcommand{\DefinitionCite}{Def.\,}

\newcommand{\norm}{\xi}   

\begin{document}

\title{Drinfeld modular forms of arbitrary rank\\
Part III: Examples}

\author{Dirk Basson \and Florian Breuer$^{1,2}$ \and Richard Pink$^2$}

\footnotetext[1]{Supported by the Alexander von Humboldt foundation, and by the NRF grant BS2008100900027.}
\footnotetext[2]{Supported through the program ``Research in Pairs'' by Mathematisches Forschungsinstitut Oberwolfach in 2010.}

\date{May 27, 2018}
\maketitle

\centerline{To Mira and Thomas}
\bigskip

\begin{abstract}
This is the third part of a series of articles providing a foundation for the theory of Drinfeld modular forms of arbitrary rank. In the present article we construct and study some examples of Drinfeld modular forms. In particular we define Eisenstein series, as well as the action of Hecke operators upon them, coefficient forms and discriminant forms. In the special case $A=\BF_q[t]$ we show that all modular forms for $\GL_r(\Gamma(t))$ are generated by certain weight one Eisenstein series, and all modular forms for $\GL_r(A)$ and $\SL_r(A)$ are generated by certain coefficient forms and discriminant forms. We also compute the dimensions of the spaces of such modular forms. 
\end{abstract}

{\advance\baselineskip by -6pt
\tableofcontents
}


\externaldocument{Part_1}
\externaldocument{Part_2}
\setcounter{section}{12}

\Bigskip


\newpage
\noindent {\Large\bf Introduction}
\bigskip

This is part III of a series of articles together with \cite{BBP1} and \cite{BBP2}, whose aim is to provide a foundation for the theory of Drinfeld modular forms of arbitrary rank. Part I developed the basic analytic theory, including $u$-expansions and holomorphy at infinity. In Part II we identified the analytic modular forms from Part I with the algebraic modular forms defined in \cite{Pink}. In the present Part III we illustrate the general theory by constructing some important families of modular forms. 

The sections of all parts are numbered consecutively, thus Sections 1--6 appear in Part I and Sections 7--12 in Part II. All the definitions and notation from Parts I and II remain in force, and we refer to proclamations in the other parts without any special indication.

For the convenience of the reader we recall some definitions and notations which we will need to describe the contents of this paper.

\medskip
Let $F$ be a global function field with exact field of constants $\BF_q$ of cardinality $q$. Fix a place $\infty$ of $F$ and denote by $A$ the ring of elements of $F$ which are regular away from~$\infty$. Let $\Finf$ be the completion of $F$ at $\infty$, and let $\Cinf$ denote the completion of an algebraic closure of~$\Finf$. Fix an integer $r \geq 2$ and an auxiliary constant $\xi\in\Cinf^\times$.
As in (\ref{DefOfOmega}) we identify the Drinfeld period domain of rank $r$ over~$\Cinf$ with the set of column vectors
$$\Omega^r := \bigl\{ (\omega_1,\ldots,\omega_r)^T \in \Cinf^r \bigm|
\omega_1,\ldots,\omega_r\ \Finf\hbox{-linearly independent and\ }\omega_r=\norm
\bigr\}.$$
Let $L$ be a finitely generated projective $A$-submodule of rank $r$ of~$F^r$, viewed as a set of row vectors. For any $\omega\in\Omega^r$ we thus obtain a strongly discrete $A$-lattice $L\omega\subset\Cinf$ of rank~$r$. Our convention on row vectors implies that $\GL_r(F)$ acts on $F^r$ from the right. We denote the stabiliser of $L$ by
$$\Gamma_L\ :=\ \bigl\{ \gamma\in\GL_r(F) \bigm| L\gamma = L\bigr\}.$$
For $L=A^r$ we simply have $\Gamma_L=\GL_r(A)$. Note that for any non-zero ideal $N\subset A$, an element of $\GL_r(F)$ stabilises the lattice $L$ if and only if it stabilises the lattice $N^{-1}L$; thus $\Gamma_{N^{-1}L}=\Gamma_L$. More generally, for any coset $v+L\subset F^r$ we consider the congruence subgroup 
$$\Gamma_{v+L}\ :=\ 
\bigl\{ \gamma\in\GL_r(F) \bigm| v\gamma+L\gamma = v+L\bigr\}
\ <\ \Gamma_L.$$
Also, for any non-zero ideal $N\subset A$ we consider the principal congruence group
$$\Gamma_L(N)\ :=\ \bigcap_{v\in N^{-1}L}\Gamma_{v+L}
\ =\ \ker\bigl(\Gamma_L\to\Aut(N^{-1}L/L)\bigr).$$
All these groups are arithmetic subgroups of $\GL_r(F)$.

For any $\gamma\in\GL_r(\Finf)$ and $\omega\in\Omega^r$, we let 
\[j(\gamma,\omega) := \norm^{-1}\cdot(\hbox{last entry of\ }\gamma\omega) \in \Cinf^\times. \tag{\ref{DefOfj}}\]
Then matrix multiplication from the left followed by a normalisation
defines a left action of $\GL_r(\Finf)$ on $\Omega^r$ via
\[\gamma(\omega) := j(\gamma,\omega)^{-1}\gamma\omega.\tag{\ref{DefOfActionOnOmega}}\] 
For any integers $k$ and $m$, the group $\GL_r(\Finf)$ acts on the space of functions $f : \Omega^r \to \Cinf$ via the Petersson slash operator, as follows:
\[(f|_{k,m}\gamma)(\omega)
:= \det(\gamma)^{m}j(\gamma,\omega)^{-k}f\big(\gamma(\omega)\big). \tag{\ref{DefOfActionOnF}}\]
For any arithmetic subgroup $\Gamma < \GL_r(F)$ we defined the space of weak modular forms of weight $k$ and type $m$ for $\Gamma$ in Definition \ref{Def:WeakModForm} as
\[\CW_{k,m}(\Gamma) := \{ f : \Omega^r \to \Cinf \;|\; \text{$f$ is holomorphic and $f|_{k,m} = f$}\}.\]
When $m=0$ or $\det(\Gamma)=\{1\}$ we suppress the $m$ from our notation. 

To define modular forms, we constructed $u$-expansions of weak modular forms in Section \ref{Sec:HolomorphyB} as follows. 
Write an element of $\Omega^r$ in the form $\omega = \binom{\omega_1}{\omega'}$ with $\omega_1\in\Cinf$ and $\omega'\in\Omega^{r-1}$. Set
\[u_{\omega'}(\omega_1) := \frac{1}{e_{\Lambda'\omega'}(\omega_1)} \in\Cinf^\times, \tag{\ref{DefOfuParameter}}\]
where $e_{\Lambda'\omega'}$ is the exponential function associated to the strongly discrete subgroup $\Lambda'\omega'\subset\Cinf$ and $\Lambda'\subset F^{r-1}$ is a certain group of column vectors associated to $\Gamma$ by (\ref{DefLambdaPrime}).
In Proposition \ref{ModFormsLaurentExpansion} and Theorem \ref{Thm:CoefficientsAreModularForms} we proved that every weak modular form $f : \Omega^r \to \Cinf$ of weight $k$ admits a unique $u$-expansion
\[
f(\omega) = \sum_{n\in\BZ}f_n(\omega')u_{\Lambda'\omega'}(\omega_1)^n
\]
which converges on a suitable neighbourhood of infinity, and where the coefficients $f_n$ are weak modular forms of weight $k-n$ for an associated arithmetic group $\Gamma_M < \GL_{r-1}(F)$. We say that $f$ is holomorphic at infinity if $f_n \equiv 0$ for all $n<0$, and that $f$ vanishes at infinity if $f_n \equiv 0$ for all $n\leq 0$. The spaces of modular forms, resp.\ cusp forms of weight $k$ and type $m$ for $\Gamma$ are then defined by
\begin{eqnarray*}
\CM_{k,m}(\Gamma) & := & \{ f \in \CW_{k,m}(\Gamma) \;|\; \text{$f|_{k,m}\gamma$ is holomorphic at infinity for all $\gamma\in\GL_r(F)$}\} \quad\text{and} \\
\CS_{k,m}(\Gamma) & := & \{ f \in \CW_{k,m}(\Gamma) \;|\; \text{$f|_{k,m}\gamma$ vanishes at infinity for all $\gamma\in\GL_r(F)$}\}, \quad\text{respectively.}
\end{eqnarray*}
These are finite-dimensional $\Cinf$-vector spaces, by Theorem \ref{FiniteDimension}. Again we leave out the index $m$ if $m=0$. 

\subsubsection*{Outline  of this paper}

In Section \ref{sec:Eisenstein} we construct the Eisenstein series of all weights $k\geq 1$ associated to all cosets $v+L$ and compute their $u$-expansions in Proposition \ref{EisensteinUExpansion}. In Theorem \ref{EisensteinModularForm} we show that they are modular forms of weight $k$ for the groups $\Gamma_{v+L}$.

In Section \ref{HeckeEis} we determine the action of Hecke operators (defined in Section \ref{HeckeOps}) on Eisenstein series, restricting ourselves to Hecke operators that are supported away from the level of the Eisenstein series (see Assumption \ref{EisenHeckeActionAss}). In each case, Theorem \ref{EisenHeckeAction1} identifies the Hecke image of an Eisenstein series as a linear combination of Eisenstein series. In particular, we deduce that Eisenstein series are eigenforms under many Hecke operators.

Coefficient forms are defined in Section \ref{CoeffForms}, they are modular forms for $\Gamma_L$ which occur as coefficients of Drinfeld modules, isogenies or exponential functions associated to the lattice $L\omega$. 

Section \ref{DiscForms} deals with discriminant forms, which arise as highest coefficients of Drinfeld modules or as roots thereof. These are always cusp forms. Certain $(q-1)$-st roots are examples of modular forms with non-zero type $m$. 

Lastly, we discuss the special case of $A=\BF_q[t]$ and $L=A^r$ in Section \ref{SecFqt}. Here we exploit the explicit description of algebraic modular forms for $\Gamma(t)$ from \cite{PinkSchieder} and \cite{Pink} together with our identification of analytic and algebraic modular forms from Part~II. This allows us to prove in Theorem \ref{ModRing1} that the graded ring $\CM_*(\Gamma(t))$ of modular forms of all weights for $\Gamma(t)$ is generated over $\Cinf$ by the weight one Eisenstein series $E_{1,v+L}$ for all $v\in t^{-1}L\smallsetminus L$. Using invariants, we then deduce that the rings $\CM_*(\GL_r(A))$ and $\CM_*(\SL_r(A))$ are generated by suitable algebraically independent coefficient forms. This generalises known results from the $r=2$ case due to Cornelissen, Goss and Gekeler, respectively. Lastly, we give some dimension formulae in Theorem \ref{DimensionFormulas}.

\section{Eisenstein series}
\label{sec:Eisenstein}

For any integer $k\ge1$ and any vector $v\in F^r$ we define the \emph{Eisenstein series of weight $k$ associated to the coset $v+L$} by
\UseTheoremCounterForNextEquation
\begin{equation}\label{EisensteinSeriesDef}
E_{k,v+L}(\omega)\ :=\ 
\sum_{0\,\not=\,x\,\in\,v+L} (x\omega)^{-k}.
\end{equation}

\begin{Prop}\label{EisensteinSeriesHolomorphic}
This series defines a holomorphic function $\Omega^r \to \Cinf$.
\end{Prop}

\begin{Proof}
By Proposition \ref{AffinoidCovering} it suffices to show that the series converges uniformly on the affinoid set $\Omega^r_n$ from (\ref{Omega_n}) for every~$n$. For this observe that any $x\in F^r\setminus\{0\}$ determines a unimodular $F_\infty$-linear form $\frac{x}{|x|}$ on~$\Finf^r$. For any $\omega\in\Omega_n^r$ it follows that
$$|x\omega|
\ =\ |x|\cdot\bigl|\tfrac{x}{|x|}\omega\bigr|
\ \smash{\stackrel{\eqref{DOmega}}{\geq}}\ 
|x|\cdot h(\omega)\cdot|\omega|
\ \smash{\stackrel{\eqref{Omega_n}}{\geq}}\ 
|x|\cdot|\pi^n|\cdot|\omega|
\ \smash{\stackrel{\ref{Lem:Omeganbounds}}{\geq}}\ 
|x|\cdot|\pi^n|\cdot|\xi|.$$
As $x$ runs through $(v+L)\setminus\{0\}$, the norm $|x|$ goes to infinity; hence $|x\omega|^{-k}$ goes to zero uniformly over $\Omega^r_n$, as desired.
\end{Proof}

\medskip
Some basic transformation properties of Eisenstein series are:

\begin{Prop} \label{EisensteinBasicProps}
\begin{enumerate}
\item[(a)] For every $\gamma\in\GL_r(F)$ we have $E_{k,v+L}|_k\gamma = E_{k,v\gamma+L\gamma}$. 
\item[(b)] In particular $E_{k,v+L}$ is a weak modular form of weight $k$ for the group $\Gamma_{v+L}$.
\item[(c)] For any $A$-submodule of finite index $L'\subset L$ we have
$E_{k,v+L} = \sum_{v'+L'} E_{k,v'+L'}$,
where the sum extends over all $L'$-cosets $v'+L'\subset v+L$.
\end{enumerate}
\end{Prop}

\begin{Proof}
(a) results from the calculation
\begin{eqnarray*}
(E_{k,v+L}|_k\gamma)(\omega)
\!\!&\stackrel{(\ref{DefOfActionOnF})}{=}&\!\!
j(\gamma,\omega)^{-k} \cdot \kern-8pt \sum_{0\,\not=\,x\,\in\,v+L} \kern-6pt (x\cdot\gamma(\omega))^{-k} \\
\!\!&=&\!\!
\kern-8pt \sum_{0\,\not=\,x\,\in\,v+L}\kern-6pt (j(\gamma,\omega)\cdot x\cdot\gamma(\omega))^{-k} \\
\!\!&\stackrel{(\ref{DefOfActionOnOmega})}{=}&\!\!
\kern-8pt \sum_{0\,\not=\,x\,\in\,v+L} \kern-6pt (x\gamma\omega)^{-k} \\
\!\!&=&\!\!
E_{k,v\gamma+L\gamma}(\omega).
\end{eqnarray*}
(b) is a direct consequence of (a), and (c) is obvious from the definition (\ref{EisensteinSeriesDef}). 
\end{Proof}

\medskip
Our next goal is to determine the $u$-expansion of $E_{k,v+L}$, which requires some preparation. For any strongly discrete $\BF_q$-subspace $H\subset\Cinf$ consider the power series expansion of the exponential function
\UseTheoremCounterForNextEquation
\begin{equation}\label{ExpDef13}
e_H(z)\ :=\ z\cdot\!\!\prod_{h\in H\smallsetminus\{0\}}\Bigl(1-\frac{z}{h}\Bigr)
\ =\ \sum_{i=0}^\infty e_{H,q^i}z^{q^i}
\end{equation}
with $e_{H,q^i}\in\Cinf$ and $e_{H,1}=1$ that is furnished by Proposition \ref{exp1}.

\begin{Prop}\label{GossPolynomials}
\begin{enumerate}
\item[(a)] For any strongly discrete $\BF_q$-subspace $H\subset\Cinf$ we have
$$e_H(z)^{-1} \ =\ \sum_{h\in H}(z-h)^{-1}.$$
\item[(b)] For every $k\ge1$, there exists a unique so-called \emph{Goss polynomial} $G_k(X,Y_1,Y_2,\ldots)$ with coefficients in $\BF_p$ in the variables $X$ and $Y_i$ for all integers $1\le i<\log_qk$, such that for every strongly discrete $\BF_q$-subspace $H\subset\Cinf$ we have
$$G_k\bigl(e_H(z)^{-1},\;e_{H,q},\;e_{H,q^2},\;\ldots\bigr) \ =\ \sum_{h\in H}(z-h)^{-k}.$$
\item[(c)] These polynomials further satisfy:
\begin{enumerate}
\item[(i)] $G_k$ is monic of degree $k$ in~$X$ and divisible by~$X$.
\item[(ii)] $G_1=X$ and $G_k=X\bigl(G_{k-1}+\sum_{1\le i<\log_qk}Y_iG_{k-q^i}\bigr)$ for all $k>1$.
\item[(iii)] $G_{pk}=G_k^p$.
\item[(iv)] $X^2\frac{\partial}{\partial X}G_k=kG_{k+1}$.
\end{enumerate}
\end{enumerate}
\end{Prop}

\begin{Proof}
The existence of these polynomials was first obtained by Goss in \cite[\PropositionCite \ 6.6]{GossAlg}, 
but in this generality see Gekeler \cite[\TheoremCite 2.6]{GekelerZeroes}.
%
\end{Proof}

\begin{Rem}\label{GossPolynomialsDegree}
\rm We shall see in Proposition \ref{EisensteinVanishingOrder} that the vanishing order at infinity of the Eisenstein series $E_{k,v+L}$ is controlled by the vanishing order of the Goss polynomial $G_k$ at $X=0$. By part (i) of Proposition \ref{GossPolynomials} (c) this vanishing order is $\ge 1$, and part (ii) implies that it is equal to $k$ for all $k\leq q$. 
In \cite{GekelerZeroes}, Gekeler gives a formula for the order of the Goss polynomial at $X=0$ in the case $A=\BF_p[t]$ and $H=\bar\pi A$, where $p$ is prime and $\bar\pi$ is the Carlitz period. This determines the vanishing order of the Eisenstein series in the rank 2 case for $A=\BF_p[t]$.
\end{Rem}

\begin{samepage}
\begin{Cor} \label{EisensteinWeight1}
For any $v\in F^r\setminus L$ we have 
$$E_{1,v+L}(\omega)\ =\ e_{L\omega}(v\omega)^{-1}.$$
\end{Cor}

\begin{Proof}
Direct computation using the substitution $x=v-\ell$ and Proposition \ref{GossPolynomials} (a):
$$E_{1,v+L}(\omega) 
\ =\ \sum_{0\,\not=\,x\,\in\,v+L} (x\omega)^{-1}
\ =\ \sum_{\ell\,\in\,L} (v\omega-\ell\omega)^{-1}
\ =\ e_{L\omega}(v\omega)^{-1}.$$
\vskip-20pt
\end{Proof}
\end{samepage}

\medskip
Now define $A$-submodules $L'$ and $L_1$ by the commutative diagram with exact rows
\UseTheoremCounterForNextEquation
\begin{equation}\label{L'L1Def}
\xymatrix@R-24pt{
0 \ar[r] & F^{r-1} \ar[rr] && F^r \ar[rr] && F \ar[r] & 0\\
& \cup \ar@{}[rr]|{x'\mapsto(0,x')} 
&& \cup \ar@{}[rr]|{(x_1,x')\mapsto x_1} && \cup & \\
0 \ar[r] & L' \ar[rr] && L \ar[rr] && L_1 \ar[r] & 0\rlap{.}\\}
\end{equation}
Since $L$ is finitely generated projective of rank~$r$, the $A$-modules $L'$ and $L_1$ are finitely generated projective of ranks $r-1$ and $1$, respectively. 
Also fix a subgroup $\tilde L_1\subset L$ which maps isomorphically to~$L_1$, so that $L= \tilde{L}_1 \oplus (\{0\}\times L')$.
Write $v=(v_1,v')\in F^r=F\times F^{r-1}$.

\begin{Lem}\label{EisLem1}
The subgroup $\Lambda'\subset F^{r-1}$ from (\ref{DefLambdaPrime}) that corresponds to $\Gamma_{v+L}\cap U(F)$ is the finitely generated $A$-submodule of rank $r-1$
$$\Lambda'\: =\ \bigl\{ \lambda'\in F^{r-1} \bigm| (v_1+L_1)\lambda'\subset L'\bigr\}.$$
Moreover, for any $x_1\in (v_1+L_1)\setminus\{0\}$ the inclusion $x_1\Lambda'\subset L'$ has finite index. 
\end{Lem}

\begin{Proof}
For any $\lambda'\in F^{r-1}$ and $(x_1,x')\in F^r=F\times F^{r-1}$ we have ${(x_1,x') \binom{1\ \lambda'}{0\kern5pt 1\kern2pt}} = {(x_1,x_1\lambda'+x')}$. By the definition of $\Gamma_{v+L}$ in the introduction
it follows that $\lambda'\in\Lambda'$ if and only if for every $(x_1,x')\in v+L$ we have $(0,x_1\lambda')\in L$, or equivalently $x_1\lambda'\in L'$. As $(x_1,x')$ runs through $v+L$, its first component $x_1$ runs through $v_1+L_1$, so the formula for $\Lambda'$ follows.

Since $L'$ and $L_1$ are finitely generated $A$-modules of ranks $r-1$ and $1$, respectively, the formula implies that $\Lambda'$ is a finitely generated $A$-submodule of rank $r-1$. For $x_1\in (v_1+L_1)\setminus\{0\}$ it follows that  $x_1\Lambda'\subset L'$ is an inclusion of finitely generated $A$-modules of the same rank and hence of finite index.
\end{Proof}

\medskip
As before we write $\omega = \binom{\omega_1}{\omega'}\in\Omega^r \subset \Cinf\times\Omega^{r-1}$. Then the expansion parameter from (\ref{DefOfuParameter}) is the function $u := u_{\omega'}(\omega_1) := e_{\Lambda'\omega'}(\omega_1)^{-1}$.

\begin{Prop}\label{EisensteinUExpansion}
We have 
$$E_{k,v+L}({\textstyle\binom{\omega_1}{\omega'}})
\ =\ \sum_{x=(x_1,x')\,\in\,v+\tilde L_1}
\left\{
\begin{array}{ll}
E_{k,x'+L'}(\omega') 
& \hbox{if $x_1=0$,} \\[5pt]
G_k\bigl(e_{L'\omega'}(x\omega)^{-1},\; e_{L'\omega',q},\; e_{L'\omega',q^2}, \; \ldots \bigr) 
& \hbox{if $x_1\not=0$,}
\end{array}
\right.$$
where $G_k$ is the $k$-th Goss polynomial from Proposition \ref{GossPolynomials} and in the second case
$$e_{L'\omega'}(x\omega)^{-1}
\ =\ \frac{u^{[L':x_1\Lambda']}}{x_1}\cdot
\frac
{\displaystyle\prod_{{\ell'\in L'\setminus x_1\Lambda'}\ {{\rm mod}\ x_1\Lambda'}}
\kern-20pt e_{\Lambda'\omega'}(x_1^{-1}\ell'\omega')\kern30pt}
{\displaystyle\prod_{{\ell'\in L'}\ {{\rm mod}\ x_1\Lambda'}}
\kern-10pt \bigl(1-e_{\Lambda'\omega'}(x_1^{-1}(\ell'-x')\omega')\cdot u\bigr)}.$$
Moreover, the right hand side converges locally uniformly for all $(u,\omega')$ in a suitable tubular neighbourhood of $\{0\}\times\Omega^{r-1}$.
\end{Prop}

\begin{Proof}
Using the fact that $L=\tilde L_1\oplus (\{0\}\times L')$, we break up the series defining $E_{k,v+L}$ as 
\UseTheoremCounterForNextEquation
\begin{equation}\label{EisensteinUExpansion1}
E_{k,v+L}(\omega)
\ =\ \sum_{0\,\not=\,x\,\in\,v+L} (x\omega)^{-k}
\ =\ \sum_{x\,\in\,v+\tilde L_1}
\biggl(\sum_{0\,\not=\,y\,\in\,x+(\{0\}\times L')} (y\omega)^{-k}\biggr).
\end{equation}
Write $x=(x_1,x')\in F^r=F\times F^{r-1}$, and observe that for any $y=(y_1,y')\in F^r=F\times F^{r-1}$ we have $y\omega = y_1\omega_1+y'\omega'$. 

If $x_1=0$, the inner sum of (\ref{EisensteinUExpansion1}) is just
$$\sum_{0\,\not=\,y'\,\in\,x'+L'} (y'\omega')^{-k}
\ =\ E_{k,x'+L'}(\omega').$$
Such a term occurs only if $v$ lies in $L+(\{0\}\times F^{r-1})$, and then it occurs for a unique~$x$. 

If $x_1\not=0$, we write $y=x-(0,\ell')$, so that $y\omega = x\omega-\ell'\omega'$. By Proposition \ref{GossPolynomials} (b) the inner sum of (\ref{EisensteinUExpansion1}) then becomes
$$\sum_{\ell'\in L'} (x\omega-\ell'\omega')^{-k}
\ =\ G_k\bigl(e_{L'\omega'}(x\omega)^{-1},\; e_{L'\omega',q},\; e_{L'\omega',q^2}, \; \ldots \bigr).$$
To transform $e_{L'\omega'}(x\omega)$ we proceed as in the proof of Proposition \ref{GenDrinModCodim1}. First, by Lemma \ref{EisLem1} we have an inclusion of finite index $\Lambda'\omega' \subset x_1^{-1}L'\omega'$, and by the $\Finf$-linear independence of the coefficients of $\omega$ the index is precisely ${[L':x_1\Lambda']}$. By the additivity of the exponential function we have
$$e_{\Lambda'\omega'}(x_1^{-1}x\omega)
\ =\ e_{\Lambda'\omega'}(\omega_1+x_1^{-1}x'\omega')
\ =\ u^{-1}+e_{\Lambda'\omega'}(x_1^{-1}x'\omega')$$
with $u=e_{\Lambda'\omega'}(\omega_1)^{-1}$. 
Using Proposition \ref{exp2} we deduce that
\begin{eqnarray*}
e_{L'\omega'}(x\omega)
&\!\!=\!\!& x_1\cdot e_{x_1^{-1}L'\omega'}(x_1^{-1}x\omega) \\
&\!\!=\!\!& x_1\cdot
e_{e_{\Lambda'\omega'}(x_1^{-1}L'\omega')}\bigl(
e_{\Lambda'\omega'}(x_1^{-1}x\omega)\bigr) \\
&\!\!=\!\!& x_1\cdot
e_{e_{\Lambda'\omega'}(x_1^{-1}L'\omega')}\bigl(
u^{-1}+e_{\Lambda'\omega'}(x_1^{-1}x'\omega')\bigr).
\end{eqnarray*}
By the definition and the additivity of the exponential function this in turn yields
\begin{eqnarray*}
e_{L'\omega'}(x\omega)
&\!\!=\!\!& x_1\cdot
\bigl(u^{-1}+e_{\Lambda'\omega'}(x_1^{-1}x'\omega')\bigr) \cdot\kern-10pt
\prod_{{\ell'\in L'\setminus x_1\Lambda'}\atop{{\rm modulo}\ x_1\Lambda'}}\!\!
\left(1 - \frac{u^{-1}+e_{\Lambda'\omega'}(x_1^{-1}x'\omega')}
{e_{\Lambda'\omega'}(x_1^{-1}\ell'\omega')}\right) \\
&\!\!=\!\!& x_1\cdot
\bigl(u^{-1}+e_{\Lambda'\omega'}(x_1^{-1}x'\omega')\bigr) \cdot\kern-10pt
\prod_{{\ell'\in L'\setminus x_1\Lambda'}\atop{{\rm modulo}\ x_1\Lambda'}}\kern-10pt
\frac{e_{\Lambda'\omega'}(x_1^{-1}(\ell'-x')\omega')-u^{-1}}
{e_{\Lambda'\omega'}(x_1^{-1}\ell'\omega')} \\
&\!\!=\!\!& x_1\cdot
\frac{1+e_{\Lambda'\omega'}(x_1^{-1}x'\omega')\cdot u}
{u^{[L':x_1\Lambda']}}\cdot\kern-10pt
\prod_{{\ell'\in L'\setminus x_1\Lambda'}\atop{{\rm modulo}\ x_1\Lambda'}}\kern-10pt
\frac{e_{\Lambda'\omega'}(x_1^{-1}(\ell'-x')\omega')\cdot u-1}
{e_{\Lambda'\omega'}(x_1^{-1}\ell'\omega')} \\
&\!\!=\!\!& \frac{x_1}{u^{[L':x_1\Lambda']}}\cdot
\frac
{\displaystyle\prod_{{\ell'\in L'}\ {{\rm mod}\ x_1\Lambda'}}
\kern-10pt \bigl(1-e_{\Lambda'\omega'}(x_1^{-1}(\ell'-x')\omega')\cdot u\bigr)}
{\displaystyle\prod_{{\ell'\in L'\setminus x_1\Lambda'}\ {{\rm mod}\ x_1\Lambda'}}
\kern-20pt e_{\Lambda'\omega'}(x_1^{-1}\ell'\omega')\kern30pt},
\end{eqnarray*}
where the last transformation uses the fact that $(-1)^{[L':x_1\Lambda']-1}=1$ because $[L':x_1\Lambda']$ is a power of~$q$. Combining everything we obtain the desired formula.

For the convergence take any $n>0$. By Proposition \ref{NewEarth} (c) there exists a constant $c_n>0$, such that for any $\omega'\in\Omega^{r-1}_n$ and any $x'\in F_\infty^{r-1}$ we have $|e_{\Lambda'\omega'}(x'\omega')|<c_n$. In particular this inequality holds for $x_1^{-1}\ell'$ and $x_1^{-1}(\ell'-x')$ in place of~$x'$. Thus if $|u|\le r_n := (2c_n)^{-1}$, we have $|e_{\Lambda'\omega'}(x_1^{-1}(\ell'-x')\omega')\cdot u|<2^{-1}$, so the geometric series for 
$$\frac{1}{1-e_{\Lambda'\omega'}(x_1^{-1}(\ell'-x')\omega')\cdot u}$$
converges uniformly to a value of norm~$1$. Combining the inequalities yields the bound
$$\left|\frac{u^{[L':x_1\Lambda']}}{x_1}\cdot
\frac
{\displaystyle\prod_{{\ell'\in L'\setminus x_1\Lambda'}\ {{\rm mod}\ x_1\Lambda'}}
\kern-20pt e_{\Lambda'\omega'}(x_1^{-1}\ell'\omega')\kern30pt}
{\displaystyle\prod_{{\ell'\in L'}\ {{\rm mod}\ x_1\Lambda'}}
\kern-10pt \bigl(1-e_{\Lambda'\omega'}(x_1^{-1}(\ell'-x')\omega')\cdot u\bigr)}
\right|
\ \le\ \frac{r_n^{[L':\ell_1\Lambda']}c_n^{[L':\ell_1\Lambda']-1}}{|x_1|}
\ =\ \frac{2^{-[L':\ell_1\Lambda']}}{|x_1|c_n}.$$
Also recall that $G_k$ is a polynomial of fixed degree in~$X$ which is divisible by~$X$, and the values $e_{L'\omega',q},\; e_{L'\omega',q^2}, \; \ldots$ for the other variables are holomorphic functions on $\Omega^{r-1}$ and hence bounded on $\Omega^{r-1}_n$. As both $|x_1|$ and $[L':x_1\Lambda']$ go to infinity with~$x_1$, this proves that the right hand side of the formula for $E_{k,x'+L'}(\omega')$ converges uniformly for all $(u,\omega') \in B(0,r_n)\times\Omega^{r-1}_n$. Varying $n$ it therefore converges locally uniformly on the tubular neighbourhood $\bigcup_{n\ge1}B(0,r_n)\times\Omega^{r-1}_n$.
\end{Proof}

\begin{Rem}\label{EisensteinUExpansionRem}
\rm In principle, the $u$-expansion of $E_{k,v+L}$ in terms of powers of $u$ can be computed from Proposition \ref{EisensteinUExpansion} by multiplying out the geometric series involved. As it stands, however, the sum is essentially a sum over a coset of $L_1\subset F$, which is a fractional ideal of~$A$. In the rank 2 case, Petrov \cite{Petrov} has shown that there are many Drinfeld modular forms with such expansions and that they exhibit many desirable properties because of it. One may ask if there are other examples in the higher rank case.
\end{Rem}


\begin{Prop}\label{EisensteinVanishingOrder}
\begin{enumerate}
\item[(a)] The $u$-expansion of $E_{k,v+L}(\omega)$ has constant term $E_{k,x'+L'}(\omega')$ if $v\in L+(0,x')$ for some $x'\in F^{r-1}$, and constant term $0$ otherwise.
\item[(b)] If $v\not\in L+(\{0\}\times F^{r-1})$, the order at infinity of $E_{k,v+L}(\omega)$ with respect to the group $\Gamma_{v+L}\cap U(F)$ is at least 
$$\ord_X(G_k)\cdot \min\bigl\{ [L':x_1\Lambda'] \bigm| x_1\in v_1+L_1 \bigr\}.$$
\end{enumerate}
\end{Prop}

\begin{Proof}
Assertion (a) follows from Proposition \ref{EisensteinUExpansion} and the fact that the Goss polynomial $G_k$ is divisible by~$X$. In (b) let $d := \ord_X(G_k)$ denote the vanishing order at $X=0$ of $G_k$ as a polynomial in independent variables $X,Y_1,Y_2,\ldots$ and write $G_k=X^dH(Y_1,Y_2,\ldots)+($higher terms in~$X$). Then each summand in Proposition \ref{EisensteinUExpansion} contributes 
$$\left(
u^{[L':x_1\Lambda']}\cdot
\frac{\displaystyle\prod_{{\ell'\in L'\setminus x_1\Lambda'}\ {{\rm mod}\ x_1\Lambda'}}
\kern-20pt e_{\Lambda'\omega'}(x_1^{-1}\ell'\omega')}{x_1}
\right)^d
\cdot H\bigl(e_{L'\omega',q},\; e_{L'\omega',q^2},\ldots\bigr)
+\bigl(\hbox{higher terms in $u$}\bigr)$$
to the $u$-expansion of $E_{k,v+L}(\omega)$. Recall that $v=(v_1,v')$, so that as $x=(x_1,x')$ runs through $v+\tilde L_1$, its first component $x_1$ runs through $v_1+L_1$. Combining this yields the desired lower bound.
\end{Proof}

\begin{Rem}\label{EisensteinVanishingOrderRem}
\rm For the purposes explained in Remark \ref{FundDiscRem} below, one should hope that the inequality in Proposition \ref{EisensteinVanishingOrder} is always an equality in the case $k=1$. By (\ref{DiscFormDef2}) this would yield a formula for the order at infinity of every discriminant form. For example we have:
\end{Rem}

\begin{Prop}\label{EisensteinVanishingOrderEx}
If $A=\BF_q[t]$, for any $v\in t^{-1}L\setminus L$ the order at infinity of $E_{1,v+L}$ with respect to the group $\Gamma_{v+L}\cap U(F)$ is $0$ if $v\in L+(\{0\}\times F^{r-1})$ and $1$ otherwise.
\end{Prop}

\begin{Proof}
As above write $v=(v_1,v')$. If $v_1\in L_1$, the $u$-expansion of $E_{1,v+L}$ has constant term $E_{1,v'+L'}$ by Proposition \ref{EisensteinVanishingOrder} (a), which is non-zero by Corollary \ref{EisensteinWeight1}; hence the order is $0$ in this case. 

Otherwise we have $t^{-1}L_1=\BF_q\cdot v_1+L_1$ and this $A$-module is generated by a unique element $x_1\in v_1+L_1$. By Lemma \ref{EisLem1} we deduce that $\Lambda' = x_1^{-1}L'$. This $x_1$ is then the unique element of the coset $v_1+L_1$ that satisfies $[L':x_1\Lambda']=1$. Since, moreover, $G_1(X)=X$ by Proposition \ref{GossPolynomials} (b), Proposition \ref{EisensteinUExpansion} implies that $E_{1,v+L}(\omega) = \frac{u}{x_1}+($higher terms in~$u$). The order is therefore $1$ in that case.
\end{Proof}

\begin{Thm}\label{EisensteinModularForm}
The Eisenstein series $E_{k,v+L}$ is a modular form of weight $k$ for the group $\Gamma_{v+L}$.
\end{Thm}

\begin{Proof}
By Proposition \ref{EisensteinBasicProps} (b) it is already a weak modular form for $\Gamma_{v+L}$. Moreover, for every $\gamma\in\GL_r(F)$ we have $E_{k,v+L}|_k\gamma = E_{k,v\gamma+L\gamma}$ by Proposition \ref{EisensteinBasicProps} (a), and the latter is holomorphic at infinity by Proposition \ref{EisensteinUExpansion}.
\end{Proof}

\section{Hecke action on Eisenstein series}
\label{HeckeEis}

For any coset $v+L$ the quotient $(Av+L)/L$ is a finite $A$-module that is generated by one element; hence it is isomorphic to $A/N$ for a unique non-zero ideal~$N$. Equivalently $N$ is the largest ideal of $A$ such that $\Gamma_{v+L}$ contains the principal congruence subgroup $\Gamma_L(N)$. We can therefore view $N$ as a kind of \emph{level} of the Eisenstein series $E_{k,v+L}$. In this section we compute the effect on $E_{k,v+L}$ of a Hecke operator that is supported away from~$N$.

\medskip
For any finitely generated $A$-submodule $L\subset F^r$ of rank $r$ and any prime $\Fp\subset A$ let $L_\Fp$ denote the closure of $L$ in~$F_\Fp^r$, which is a finitely generated $A_\Fp$-submodule of rank~$r$. Note that $L$ can be recovered from the submodules $L_\Fp$ for all $\Fp$ as the intersection $F^r\cap\prod_\Fp L_\Fp$ within $(\BAfF)_{\vphantom{t}}^r$. 
Consider finitely generated projective $A$-submodules $L$, $L'\subset F^r$ of rank~$r$, vectors $v$, $v'\in F^r$, and an element $\delta\in\GL_r(F)$, which together satisfy:

\begin{Ass}\label{EisenHeckeActionAss}
For every prime $\Fp\subset A$ we have:
\begin{enumerate}
\item[(a)] $v\delta+L_\Fp\delta \subset v'+L'_\Fp$, 
\item[(b)] $v\delta+L_\Fp\delta = v'+L'_\Fp$ whenever $v\not\in L_\Fp$, and
\item[(c)] $L_\Fp\delta\not\subset\Fp L'_\Fp$.
\end{enumerate}
\end{Ass}

Here (a) is equivalent to $v\delta+L\delta\subset v'+L'$, which includes the fact that $L\delta\subset L'$. Given (a), condition (b) means that $E_{k,v+L}$ and $E_{k,v'+L'}$ are Eisenstein series of the same level~$N$ and that $T_\delta$ is supported only at primes not dividing~$N$. Property (c) is equivalent to $L\delta\not\subset\Fp L'$ for any prime~$\Fp$, which serves as normalisation. 
If $L=L'=A^r$, then (a) means that $\delta$ has coefficients in $A$ and maps $v$ into $v'+A^r$. Then, in addition, condition (b) means that the determinant of $\delta$ is relatively prime to~$N$, and (c) means that $\delta$ is not congruent to the zero matrix modulo any prime of~$A$. 
Assumption \ref{EisenHeckeActionAss} will remain in force until Theorem \ref{EisenHeckeAction1} below.


\medskip
To begin with we abbreviate
$$\begin{array}{ll}
\Gamma'\! &:=\ \Gamma_{v'+L'}, \\[3pt]
\Gamma  &:=\ \delta^{-1}\Gamma_{v+L}\delta\cap\Gamma_{v'+L'}
\ =\ \Gamma_{v\delta+L\delta}\cap\Gamma_{v'+L'}
\ <\ \Gamma'.
\end{array}$$
For any prime $\Fp\subset A$ we consider the open compact subgroups
$$\begin{array}{rl}
K'_\Fp &:=\ 
\bigl\{\; k\in\GL_r(F_\Fp) \bigm| v'k+L'_\Fp k = v'+L'_\Fp\;\bigr\}, \\[3pt]
K_\Fp &:=\ 
\bigl\{\; k\in\GL_r(F_\Fp) \bigm| v'k+L'_\Fp k = v'+L'_\Fp
\hbox{\ and\ } v\delta k+L_\Fp\delta k = v\delta+L_\Fp\delta\;\bigr\}
\ <\ K'_\Fp.
\end{array}$$
Since $L'/L\delta$ is finite, for any prime $\Fp$ not dividing its annihilator we have $L_\Fp\delta = L'_\Fp$ and hence $v\delta+L_\Fp\delta = v'+L'_\Fp$. Thus for almost all $\Fp$ we have $K_\Fp=K'_\Fp$. By Assumption \ref{EisenHeckeActionAss} (b) this is so in particular if $v\not\in L_\Fp$. Also, the equalities $L'=F^r\cap\prod_\Fp L'_\Fp$ and $L=F^r\cap\prod_\Fp L_\Fp$ imply that $\Gamma'=\GL_r(F)\cap \prod_\Fp K'_\Fp$ and $\Gamma=\GL_r(F)\cap\prod_\Fp K_\Fp$.

\begin{Lem}\label{EisenHeckeLem1}
For every $\Fp$ we have $\det(K_\Fp)=\det(K'_\Fp)$.
\end{Lem}

\begin{Proof}
If $v\not\in L_\Fp$, this follows from the fact that $K'_\Fp=K_\Fp$. Otherwise by assumption we have $L_\Fp\delta=v\delta+L_\Fp\delta \subset v'+L'_\Fp=L'_\Fp$ and both are free $A_\Fp$-modules of rank $r$ within~$F_\Fp^r$. To prove the desired statement we can conjugate everything by an arbitrary element of $\GL_r(F)$. By the elementary divisor theorem we may thus without loss of generality assume that $L'_\Fp=A_\Fp^r$ and that $L_\Fp\delta = A_\Fp^r h$ for some diagonal matrix~$h\in\GL_r(F_\Fp)$.
For any $a\in A_\Fp^\times$ the diagonal matrix $\diag(1,\ldots,1,a)$ then lies in $K_\Fp$ with determinant~$a$; hence $A_\Fp^\times < \det(K_\Fp)$. As $A_\Fp^\times$ is the unique largest compact subgroup of $F_\Fp^\times$, it follows that $\det(K_\Fp)=\det(K'_\Fp)=A_\Fp^\times$, as desired.
\end{Proof}

\begin{Lem}\label{EisenHeckeLem2}
There is a natural bijection
$$\xymatrix@R-22pt@C+30pt{
\Gamma\backslash\Gamma' \ar[r] & \prod_\Fp K_\Fp\backslash K'_\Fp, \\
\Gamma\gamma \ar@{|->}[r] & (K_\Fp\gamma)_\Fp.\\}$$
\end{Lem}

\begin{Proof}
If two cosets $\Gamma\gamma_1$ and $\Gamma\gamma_2$ have the same image, we have $K_\Fp\gamma_1=K_\Fp\gamma_2$ and hence $\gamma_1\gamma_2^{-1}\in K_\Fp$ for all~$\Fp$. Thus $\gamma_1\gamma_2^{-1}\in \GL_r(F)\cap \prod_\Fp K_\Fp = \Gamma$, and so $\Gamma\gamma_1=\Gamma\gamma_2$. The map is therefore injective. 
For the surjectivity consider any collection of cosets ${K_\Fp k_\Fp\subset K'_\Fp}$. By 
Lemma \ref{EisenHeckeLem1} we may without loss of generality assume that $k_\Fp\in \SL_r(F_\Fp)\cap K'_\Fp$. By strong approximation in the group $\SL_r$ there then exists an element $\gamma\in\SL_r(F)\cap\prod_\Fp K_\Fp k_\Fp$. This element lies in $\GL_r(F)\cap\prod_\Fp K'_\Fp = \Gamma'$; hence the map is surjective.
\end{Proof}

\medskip
Next observe that for any $\gamma\in\Gamma'$ the subset $v\delta\gamma+L\delta\gamma\subset F^r$ depends only on the coset $\Gamma\gamma$. For any $x\in F^r$ we let $C(x)$ denote the number of such cosets for which $x\in v\delta\gamma+L\delta\gamma$. 
Similarly, for any $k\in K'_\Fp$ the subset $v\delta k+L_\Fp\delta k\subset F_\Fp^r$ depends only on the coset $K_\Fp k$. For any $x\in F_\Fp^r$ we let $C_\Fp(x)$ denote the number of such cosets for which $x\in v\delta k+L_\Fp\delta k$.
For any fixed $x\in F^r$ the module $(Ax+Av'+L')/L\delta$ is finite, so for any prime $\Fp$ not dividing its annihilator we have $x\in v'+L'_\Fp = v\delta+L_\Fp\delta$ and $K'_\Fp=K_\Fp$ and hence $C_\Fp(x)=1$.

\begin{Lem}\label{EisenHeckeLem3}
For any $x\in F^r$ we have $C(x) = \prod_\Fp C_\Fp(x)$.
\end{Lem}

\begin{Proof}
Since $v\in F^r$ and $L=F^r\cap\prod_\Fp L_\Fp$, for any $\gamma\in\Gamma'$ we have the equality $v\delta\gamma+L\delta\gamma = {F^r\cap (v\delta\gamma+\prod_\Fp L_\Fp\delta\gamma)}$ within $(\BAfF)_{\vphantom{t}}^r$. Since $x\in F^r$, it follows that $x\in v\delta\gamma+L\delta\gamma$ if and only if $x\in v\delta\gamma+L_\Fp\delta\gamma$ for all~$\Fp$. But the latter condition depends only on the coset $K_\Fp\gamma$, so the product formula follows from Lemma \ref{EisenHeckeLem2}.
\end{Proof}


\medskip
Now let $q_\Fp$ denote the order of the residue field $k(\Fp):=A/\Fp$. In principle one can give an explicit formula for $C_\Fp(x)$ as a polynomial in $q_\Fp$ with coefficients in~$\BZ$. But we are only interested in $C_\Fp(x)$ modulo~$(p)$, so we restrict ourselves to determining this residue class. Let $\charact_X$ denote the characteristic function of a subset $X\subset F_\Fp^r$.

\begin{Lem}\label{EisenHeckeLem4}
For any prime $\Fp$ consider the unique integers $\mu_{\Fp,1}\ge\ldots\ge\mu_{\Fp,r}\ge0$ such that $L'_\Fp/L_\Fp\delta \cong \bigoplus_{j=1}^r A/\Fp^{\mu_{\Fp,j}}$. Then for any $x\in F_\Fp^r$ we have 
$$C_\Fp(x) \ \equiv\ \left\{
\begin{array}{ll}
\charact_{v'+L'_\Fp}(x) 
& \hbox{if $\mu_{\Fp,1}\le1$} \\[3pt]
\charact_{L'_\Fp\setminus\Fp L'_\Fp}(x)
& \hbox{if $2\le\mu_{\Fp,1}\le\mu_{\Fp,r-1}+1$} \\[3pt]
0 & \hbox{if $\mu_{\Fp,1}\ge\mu_{\Fp,r-1}+2$}
\end{array}\right\}\ \mod\ (q_\Fp).$$
\end{Lem}

\begin{Proof}
By Assumption \ref{EisenHeckeActionAss} (a) we have $v\delta+L_\Fp\delta \subset v'+L'_\Fp$, so for any $k\in K'_\Fp$ we also have $v\delta k+L_\Fp\delta k\subset v'+L'_\Fp$. Thus $C_\Fp(x)=0$ if $x\not\in v'+L'_\Fp$. So till the end of the proof we assume that $x\in v'+L'_\Fp$.
If in addition $v\delta+L_\Fp\delta = v'+L'_\Fp$, we have $K_\Fp = K'_\Fp$ and $v\delta k+L_\Fp\delta k= v'+L'_\Fp$ and hence
$C_\Fp(x)=1$.
Till the end of the proof we therefore assume that $v\delta+L_\Fp\delta \not= v'+L'_\Fp$. By Assumption \ref{EisenHeckeActionAss} (b) this implies that 
$v\in L_\Fp$ and hence $L_\Fp\delta=v\delta+L_\Fp\delta \subsetneqq v'+L'_\Fp=L'_\Fp$. 
For ease of notation we abbreviate the chosen exponents to $\mu_i := \mu_{\Fp,i}$.
Then
$\mu_1\ge1$, and Assumption \ref{EisenHeckeActionAss} (c) requires that $\mu_r=0$.

Both $L_\Fp\delta\subset L'_\Fp$ are free $A_\Fp$-modules of rank $r$ within~$F_\Fp^r$. To prove the desired statement we can conjugate everything by an arbitrary element of $\GL_r(F)$. By the elementary divisor theorem we may thus without loss of generality assume that $L'_\Fp=A_\Fp^r$ and $L_\Fp\delta = \bigoplus_{j=1}^r \Fp^{\mu_j}\!A_\Fp$. Then $K'_\Fp = \GL_r(A_\Fp)$ and 
$$K_\Fp \ =\ h^{-1}\GL_r(A_\Fp)h \cap \GL_r(A_\Fp)
\ =\ \bigl\{\; (a_{ij})_{ij}\in \GL_r(A_\Fp)\bigm|
\forall i\ge j\colon\ a_{ij}\in\Fp^{\mu_j-\mu_i}\!A_\Fp \;\bigr\}.$$
Next observe that $\Fp^{\mu_1}L'_\Fp \subset L_\Fp\delta\subset L'_\Fp$. Consider the factor module $\bar L' := L'_\Fp/\Fp^{\mu_1}L'_\Fp = (A/\Fp^{\mu_1})^r$ and its submodule $\bar L := L_\Fp\delta/\Fp^{\mu_1}L'_\Fp = \bigoplus_{j=1}^r \Fp^{\mu_j}\!A_\Fp/\Fp^{\mu_1}\!A_\Fp$.
Then $K'_\Fp$ surjects to $\bar K' := \GL_r(A/\Fp^{\mu_1})$, and the image $\bar K<\bar K'$ of $K_\Fp<K'_\Fp$ is the stabiliser of~$\bar L$. In particular we have $[\bar K':\bar K]=[K'_\Fp:K_\Fp]$. To compute this number note that the image of $K_\Fp$ in $\GL_r(k(\Fp))$ is the parabolic subgroup 
$$P(k(\Fp)) \ :=\ \bigl\{\; (\bar a_{ij})_{ij}\in \GL_r(k(\Fp))\bigm|
\forall i\ge j\colon\ \mu_j>\mu_i \Rightarrow \bar a_{ij}=0 \;\bigr\},$$
and a straightforward calculation shows that $[\GL_r(k(\Fp)):P(k(\Fp))] \equiv 1$ modulo $(q_\Fp)$. From this we deduce that 
\UseTheoremCounterForNextEquation
\begin{equation}\label{EisenHeckeLem4For1}
[K'_\Fp:K_\Fp]\ \in\ \prod_{i\ge j} q_\Fp^{\max\{0,\mu_j-\mu_i-1\}} \cdot (1 + q_\Fp\BZ).
\end{equation}
Also, let $\bar x\in\bar L'$ denote the image of $x\in v'+L'_\Fp=L'_\Fp$. Then
$$C_\Fp(x)\ =\ 
\frac{|\{\bar k\in\bar K'\mid \bar x\in\bar L\bar k\}|}{|\bar K|}
\ =\ [K'_\Fp:K_\Fp]\cdot
\frac{|\{\bar k\in\bar K'\mid \bar x\bar k^{-1}\in\bar L\}|}{|\bar K'|}.$$
If $\bar x=0$, we deduce that $C_\Fp(x) = [K'_\Fp:K_\Fp]$. Otherwise $\bar x$ lies in the subset $\bar S_\nu := \Fp^\nu \bar L' \setminus \Fp^{\nu+1}\bar L'$ for a unique exponent $0\le\nu<\mu_1$. Since $\bar S_\nu$ is an orbit under~$\bar K'$, the last fraction is equal to the proportional size of $\bar L\cap\bar S_\nu$ versus $\bar S_\nu$; hence
\UseTheoremCounterForNextEquation
\begin{equation}\label{EisenHeckeLem4For2}
C_\Fp(x)\ =\ [K'_\Fp:K_\Fp] \cdot \frac{|\bar L\cap\bar S_\nu|}{|\bar S_\nu|}.
\end{equation}
To compute these cardinalities observe that 
$$\bar L\cap \Fp^\nu \bar L'
\ =\ \bigoplus_{j=1}^r\; (\Fp^{\mu_j\!}A_\Fp \cap\Fp^\nu\!A_\Fp)/\Fp^{\mu_1}\!A_\Fp 
\ =\ \bigoplus_{j=1}^r\; \Fp^{\max\{\mu_j,\nu\}}\!A_\Fp/\Fp^{\mu_j}\!A_\Fp$$
and hence 
$$|\bar L\cap \Fp^\nu \bar L'|
\ =\ \prod_{j=1}^r q_\Fp^{\mu_1-\max\{\mu_j,\nu\}}.$$
The same calculation with $\nu+1$ in place of $\nu$ shows that 
$$|\bar L\cap \Fp^{\nu+1} \bar L'|
\ =\ \prod_{j=1}^r q_\Fp^{\mu_1-\max\{\mu_j,\nu+1\}}.$$
Together this implies that 
$$|\bar L\cap \bar S_\nu|
\ =\  q_\Fp^{\sum_{j=1}^r(\mu_1-\max\{\mu_j,\nu\})}
- q_\Fp^{\sum_{j=1}^r(\mu_1-\max\{\mu_j,\nu+1\})}.$$
Since $\mu_1>\mu_r=0$, we certainly have $\mu_1-\max\{\mu_r,\nu\} > \mu_1-\max\{\mu_r,\nu+1\}$, so the first exponent is greater than the second. Therefore
\UseTheoremCounterForNextEquation
\begin{equation}\label{EisenHeckeLem4For3}
|\bar L\cap \bar S_\nu|
\ \in\  q_\Fp^{\sum_{j=1}^r(\mu_1-\max\{\mu_j,\nu+1\})}
\cdot(-1+q_\Fp\BZ).
\end{equation}
A similar, but simpler, computation shows that 
\UseTheoremCounterForNextEquation
\begin{equation}\label{EisenHeckeLem4For4}
|\bar S_\nu| \ \in\ q_\Fp^{r(\mu_1-\nu-1)}\cdot(-1+q_\Fp\BZ).
\end{equation}
Combining the formulas (\ref{EisenHeckeLem4For1}) through  (\ref{EisenHeckeLem4For4}) we deduce that 
$$C_\Fp(x)\ \in\ q_\Fp^{c(\nu)}\cdot (1 + q_\Fp\BZ)$$
for
$$c(\nu)\ :=\ 
\sum_{i\ge j} \max\{0,\mu_j-\mu_i-1\}
+ \sum_{j=1}^r(\mu_1-\max\{\mu_j,\nu+1\})
- r(\mu_1-\nu-1).$$
By (\ref{EisenHeckeLem4For1}), the same formula is true in the case $\bar x=0$ if we set $\nu:=\mu_1$.

It remains to find out when this exponent is greater than~$0$. Combining the terms for $i=r$ with the rest of the formula and using the fact that $\mu_r=0$ yields
$$\begin{array}{rl}
c(\nu)\ =& 
\sum_{r>i\ge j} \max\{0,\mu_j-\mu_i-1\}
+ \sum_{j=1}^r\bigl(\max\{0,\mu_j-1\}-\max\{0,\mu_j-\nu-1\}\bigr) \\[3pt]
=& 
\sum_{r>i\ge j} \max\{0,\mu_j-\mu_i-1\}
+ \sum_{j=1}^r\max\{0,\min\{\mu_j-1,\nu\}\}.
\end{array}$$
Here all summands are $\ge0$. Since $\mu_1\ge\ldots\ge \mu_r$, the first sum contains a positive term if and only if $\mu_1-\mu_{r-1}-1\ge1$, and the second sum contains a positive term if and only if $\min\{\mu_1-1,\nu\}\ge1$. 
Thus
$$\begin{array}{ll}
c(\nu)>0 & \hbox{if $\mu_1\ge\mu_{r-1}+2$ or ($\mu_1\ge2$ and $\nu\ge1$),} \\[3pt]
c(\nu)=0 & \hbox{if $\mu_1\le\mu_{r-1}+1$ and ($\mu_1\le1$ or $\nu=0$).}
\end{array}$$
Combining all the cases we conclude that 
$$\begin{array}{ll}
C_\Fp(x) = 0 & \hbox{if $x\not\in v'+L'_\Fp$,} \\[3pt]
C_\Fp(x) = 1 & \hbox{if $x\in v'+L'_\Fp = v\delta+L_\Fp\delta$,} \\[3pt]
C_\Fp(x) \equiv 0\ \mod\ (q_\Fp) & \hbox{if $x\in v'+L'_\Fp \not= v\delta+L_\Fp\delta$ and ($\mu_1\ge\mu_{r-1}+2$ or ($\mu_1\ge2$ and $x\in\Fp L'_\Fp$)),} \\[3pt]
C_\Fp(x) \equiv 1\ \mod\ (q_\Fp) & \hbox{if $x\in v'+L'_\Fp \not= v\delta+L_\Fp\delta$ and $\mu_1\le\mu_{r-1}+1$ and ($\mu_1\le1$ or $x\not\in\Fp L'_\Fp$),}
\end{array}$$
Since $v'+L'_\Fp = v\delta+L_\Fp\delta$ if and only if $L'_\Fp = L_\Fp\delta$ if and only if $\mu_1=0$, the desired formula follows.
\end{Proof}

\medskip
Now recall from Definition \ref{AnalyticHecke} that the Hecke operator associated to the double coset $\Gamma_{v+L}\delta\Gamma_{v'+L'}$ is defined by 
\UseTheoremCounterForNextEquation
\begin{equation}\label{AnalyticHeckeRepeat}
T_\delta:\ \CM_k(\Gamma_{v+L}) \longto \CM_k(\Gamma_{v'+L'}),\ 
f\longmapsto \sum\nolimits_{\gamma} f|_k\,\gamma,
\end{equation}
where $\gamma$ runs through a set of representatives of the quotient $\Gamma_{v+L}\backslash\Gamma_{v+L}\delta\Gamma_{v'+L'}$.

\begin{Thm}\label{EisenHeckeAction1}
Under Assumption \ref{EisenHeckeActionAss} consider the integers $\mu_{\Fp,i}$ from Lemma \ref{EisenHeckeLem4}. 
If $\mu_{\Fp,1}\ge\mu_{\Fp,r-1}+2$ for some $\Fp$, we have
$$T_\delta E_{k,v+L}\ =\ 
\rlap{$0$.}\phantom{\sum_{I\subset S} (-1)^{|I|}\cdot E_{k,v''+L'_I}}.$$
Otherwise let $S$ be the finite set of primes $\Fp$ for which $2\le\mu_{\Fp,1}\le\mu_{\Fp,r-1}+1$. For any subset $I\subset S$ set $L'_I := \prod_{\Fp\in I}\Fp\cdot L'$. Then $v'+L' = v''+L'$ for some element $v''\in (v'+L')\cap\bigcap_{\Fp\in S}\Fp L'_\Fp$ and 
$$T_\delta E_{k,v+L}\ =\ \sum_{I\subset S} (-1)^{|I|}\cdot E_{k,v''+L'_I}.$$
\end{Thm}

\begin{Proof}
By the construction of $\Gamma$ and $\Gamma'$ we have $T_\delta f = \sum\nolimits_{\gamma} E_{k,v+L}|_k\,\delta\gamma,$ where $\gamma$ runs through a set of representatives $\CR$ of $\Gamma\backslash\Gamma'$. Using the transformation rule from Proposition \ref{EisensteinBasicProps} (a) and the definition (\ref{EisensteinSeriesDef}) of Eisenstein series we deduce that
$$(T_\delta E_{k,v+L})(\omega)\ =\ 
\sum_{\gamma\in\CR} E_{k,v\delta\gamma+L\delta\gamma}(\omega)
\ =\ \sum_{\gamma\in\CR}\; \sum_{0\,\not=\,x\,\in\,v\delta\gamma+L\delta\gamma} \kern-10pt(x\omega)^{-k}
\ = \sum_{0\,\not=\,x\,\in\, F^r} \kern-5pt C(x)\cdot(x\omega)^{-k}.$$
Here $C(x)$ is determined by Lemmas \ref{EisenHeckeLem3} and \ref{EisenHeckeLem4}: 
If $\mu_{\Fp,1}\ge\mu_{\Fp,r-1}+2$ for some $\Fp$, we have $C(x)=0$ for all $x\in F^r$. 
Otherwise, for any prime $\Fp$ in the indicated set~$S$, we have $v\in L_\Fp$ and hence $v'\in L'_\Fp$ by Assumption \ref{EisenHeckeActionAss} (b). Thus $\Fp$ does not divide the annihilator $N$ of the coset $v'+L'/L'$. By the Chinese remainder theorem there therefore exists an element $a\in\bigcap_{\Fp\in S}\Fp$ with $a\equiv1$ modulo~$N$, and then $v'' := av'$ lies in $(v'+L')\cap\bigcap_{\Fp\in S}\Fp L'_\Fp$. For any subset $I\subset S$ we then have 
$$F^r\cap\prod_{{\rm all\ }\Fp}
\left\{\!\begin{array}{rl}
v'+L'_\Fp & \hbox{if $\Fp\not\in I$}, \\
\Fp L'_\Fp & \hbox{if $\Fp\in I$},
\end{array}\!\right\}
\ =\ v''+ \left(
F^r\cap\prod_{{\rm all\ }\Fp}
\left\{\!\begin{array}{rl}
L'_\Fp & \hbox{if $\Fp\not\in I$}, \\
\Fp L'_\Fp & \hbox{if $\Fp\in I$},
\end{array}\!\right\}\right)
\ =\ v''+L'_I.$$
Lemmas \ref{EisenHeckeLem3} and \ref{EisenHeckeLem4} then imply that
\begin{eqnarray*}
C(x) & \equiv &
\prod_{\Fp\notin S}\charact_{v'+L'_\Fp}(x)
\cdot\prod_{\Fp\in S}\charact_{L'_\Fp\smallsetminus\Fp L'_\Fp}(x)
\rlap{\qquad\qquad\qquad\hbox{modulo $(q)$}} \\
& = & \prod_{\Fp\notin S}\charact_{v'+L'_\Fp}(x)
\cdot\prod_{\Fp\in S}\;\bigl[\charact_{L'_\Fp}(x)-\charact_{\Fp L'_\Fp}(x)\bigr] \\
& = & \prod_{\Fp\notin S}\charact_{v'+L'_\Fp}(x)
\cdot\sum_{I\subset S} (-1)^{|I|}
\cdot\!\prod_{\Fp\in S\setminus I}\!\!\charact_{L'_\Fp}(x) 
\cdot\prod_{\Fp\in I}\charact_{\Fp L'_\Fp}(x) \\
& = & \sum_{I\subset S} (-1)^{|I|}
\cdot\prod_{\Fp\notin I}\charact_{v'+L'_\Fp}(x)
\cdot\prod_{\Fp\in I}\charact_{\Fp L'_\Fp}(x) \\
& = &
\sum_{I\subset S} (-1)^{|I|}\cdot \charact_{v''+L'_I}(x).
\end{eqnarray*}
The desired formula now follows from the definition (\ref{EisensteinSeriesDef}) of Eisenstein series.
\end{Proof}


\begin{Cor}\label{EisenHeckeAction1a}
Consider any $\delta\in\GL_r(F)$ such that for every prime $\Fp\subset A$ we have:
\begin{enumerate}
\item[(a)] $v\delta+L_\Fp\delta \subset v'+L'_\Fp$, 
\item[(b)] $v\delta+L_\Fp\delta = v'+L'_\Fp$ whenever $v\not\in L_\Fp$, and
\item[(c)] $\Fp L'_\Fp \subsetneqq L_\Fp\delta$.
\end{enumerate}
Then the Hecke operator $T_\delta$ associated to the double coset $\Gamma_{v+L}\delta\Gamma_{v'+L'}$ satisfies
$$T_\delta E_{k,v+L}\ =\ E_{k,v'+L'}.$$
\end{Cor}

\begin{Proof}
In that case Assumption \ref{EisenHeckeActionAss} holds with $\mu_{\Fp,1}\le 1$ for all~$\Fp$; hence we are in the second case of Theorem \ref{EisenHeckeAction1} with $S=\emptyset$.
\end{Proof}

\begin{Prop}\label{HeckeActionScalarTwist}
Consider any arithmetic subgroups $\Gamma$, $\Gamma'<\GL_r(F)$, any element $\delta\in\GL_r(F)$, and any scalar $a\in F^\times$. Then the Hecke operators $T_\delta$ and $T_{a^{-1}\delta}$ associated to the double cosets $\Gamma\delta\Gamma'$ and $\Gamma a^{-1}\delta\Gamma'$ satisfy
$$T_{a^{-1}\delta}\ =\ a^k\cdot T_\delta.$$
\end{Prop}

\begin{Proof}
As $\gamma$ runs through a set of representatives of $\Gamma\backslash\Gamma\delta\Gamma'$, the element $a^{-1}\gamma$ runs through a set of representatives of $\Gamma\backslash\Gamma a^{-1}\delta\Gamma'$. Since $f|_k(a^{-1}\gamma) = f|_k(a^{-1}\cdot\Id_r)|_k\gamma = a^k\cdot f|_k\gamma$ by (\ref{fCocycle}) and (\ref{IdrAction}), the formula follows from the definition of Hecke operators \ref{AnalyticHecke}.
\end{Proof}

\begin{Rem}\label{HeckeActionScalarTwistRem}
\rm Using Proposition \ref{HeckeActionScalarTwist}, one can express any Hecke operator in terms of another Hecke operator that is associated to a matrix with coefficients in~$A$. If one prefers, one can also require that the inverse matrix has coefficients in~$A$.
\end{Rem}

\begin{Rem}\label{EisenHeckeActionRem1}
\rm Combining Proposition \ref{HeckeActionScalarTwist} with Theorem \ref{EisenHeckeAction1} or Corollary \ref{EisenHeckeAction1a}, one obtains an explicit formula for $T_{a^{-1}\delta} E_{k,v+L}$ as well. 
In the special case $v'+L'=v+L$ one obtains many Hecke operators for which $E_{v+L}$ is an eigenform with eigenvalue $1$ or $a^k$.
\end{Rem}

\begin{Rem}\label{EisenHeckeActionRem2}
\rm In the case $r=2$ Theorem \ref{EisenHeckeAction1} was proved by Gekeler \cite[VIII.1]{GekelerDMC}.
%
For instance, for $L=L'=A^2$, the Hecke operator in \cite{GekelerDMC} associated to a prime element $\pi\in A$ is $T_\delta$ for the matrix $\delta = \left(\begin{smallmatrix}1 & 0 \\ 0 & \kern2pt\pi^{-1}\kern-3pt\end{smallmatrix}\right)$ and satisfies $T_{\delta} E_{k,L}=\pi^k\cdot E_{k,L}$.
%
%
%
%
\end{Rem}

\section{Coefficient forms}
\label{CoeffForms}

As before we fix a finitely generated projective $A$-submodule $L\subset F^r$ of rank~$r$. We will show that the coefficients of the exponential function $e_{L\omega}$ and of the associated Drinfeld $A$-module are modular forms for the group~$\Gamma_L$; these are the \emph{coefficient forms} in the title. We will also exhibit them as polynomials in Eisenstein series. 
The coefficients of $e_{L\omega}$ have been studied in a special case, for instance in \cite[II.2]{GekelerDMC} and \cite{GekelerPara}.

\medskip
For every $k\ge0$ we write $e_{k,L}(\omega) := e_{L\omega,q^k}$, so that $e_{L\omega}(z) = \sum_{k=0}^\infty e_{k,L}(\omega)z^{q^k}$ with $e_{0,L}=1$.
Then by  \cite[(9)]{BR09} we have
\UseTheoremCounterForNextEquation
\begin{equation}\label{ExpCoeff1a}
e_{k,L}\ =\ E_{q^k-1,L} + \sum_{j=1}^{k-1} e_{j,L} \cdot E_{q^{k-j}-1,L}^{q^j}.
\end{equation}
By direct calculation \cite[Lemma 3.4.13]{basson_thesis} this is equivalent to the more suggestive fact that $z-\sum_{i\ge1}E_{q^i-1,L}(\omega)z^{q^i}$ is the compositional inverse of~$e_{L\omega}$, in other words, that for all $\omega\in\Omega^r$ and $z\in\BC$ we have
\UseTheoremCounterForNextEquation
\begin{equation}\label{ExpCoeff1b}
e_{L\omega}\Bigl(z-\sum_{i\ge1}E_{q^i-1,L}(\omega)z^{q^i}\Bigr)\ =\ z.
\end{equation}
By induction on~$k$ the recursion formula (\ref{ExpCoeff1a}) implies that $e_{k,L}$ is a universal polynomial with coefficients in $\BF_p$ in the Eisenstein series $E_{q^i-1,L}$ for all $1\le i\le k$.

\begin{Prop}\label{ExpCoeff2}
For all $k\geq 0$ we have:
\begin{enumerate}
\item[(a)] $e_{k,L}|_{q^k-1}\gamma = e_{k,L\gamma}$ for all $\gamma\in\GL_r(F)$.
\item[(b)] $e_{k,L}$ is a modular form of weight $q^k-1$ for the group~$\Gamma_L$.
\item[(c)] The $u$-expansion of $e_{k,L}$ has constant term $e_{k,L'}$ with $L'$ as in (\ref{L'L1Def}). In particular $e_{k,L}$ is not a cusp form.
\end{enumerate}
\end{Prop}

\begin{Proof}
For any $\gamma\in\GL_r(F)$ the exponential function associated to the lattice $L\gamma(\omega)\subset\Cinf$ satisfies
$$e_{L\gamma(\omega)}\ =\ e_{j(\gamma,\omega)^{-1}L\gamma\omega}(z)\ \stackrel{\ref{exp2}}{=}\ j(\gamma,\omega)^{-1}e_{L\gamma\omega}\bigl(j(\gamma,\omega)z\bigr).$$
Comparing coefficients of $z^{q^k}$ in the respective power series expansions yields 
$$e_{L,q^k}\bigl(\gamma(\omega)\bigr)\ =\ j(\gamma,\omega)^{q^k-1} e_{L\gamma,q^k}(\omega),$$
proving (a).
Part (b) follows from Theorem \ref{EisensteinModularForm} and the formula (\ref{ExpCoeff1a}) by induction on~$k$.
To prove (c), write $\omega = \binom{\omega_1}{\omega'}$ as before. For any fixed $\omega'\in\Omega^{r-1}$, if $\omega_1$ goes to infinity, the defining formula (\ref{ExpDef13}) shows that $e_{L\omega}$ goes to $e_{L'\omega'}$ coefficientwise.
Thus $e_{k,L}$ goes to $e_{k,L'}$, and since the latter is non-zero, it follows that $e_{k,L}$ is not a cusp form.
\end{Proof}

\medskip
Next let $(\BGacomma{\Omega^r},\psi^{L})$ be the Drinfeld $A$-module of rank $r$ over $\Omega^r$ that was associated to $L$ in (\ref{PsiLDef}). Following (\ref{PsiLaDef}) and (\ref{ExpDef}) and Corollary \ref{EisensteinWeight1}, for any $a\in A\setminus\{0\}$ and any $\omega\in\Omega^r$ we then have 
\UseTheoremCounterForNextEquation
\begin{equation}\label{CoeffsDef1}
\psi^{L\omega}_a(X)
\ =\ a\cdot X\cdot\kern-10pt
\prod_{\textstyle{{v\,\in\,a^{-1}L\setminus L}\atop{{\rm modulo}\ L}}}
\kern-5pt \bigl(1-E_{1,v+L}(\omega)\cdot X\bigr).
\end{equation}
This is an $\BF_q$-linear polynomial of degree $[a^{-1}L:L]=q^{r\deg(a)}$ in~$X$. We expand it as
\UseTheoremCounterForNextEquation
\begin{equation}\label{CoeffsDef2}
\psi^{L\omega}_a(X)\ =\ \sum_{i\ge0} g^L_{a,i}(\omega)\cdot X^{q^i}
\end{equation}
with holomorphic functions $g^L_{a,i}$ on~$\Omega^r$, which are non-zero for $i=0$ and $i=r\deg(a)$ but zero whenever $i>r\deg(a)$. 
The formula (\ref{CoeffsDef1}) implies that each $g^L_{a,k}$ is a homogeneous symmetric polynomial of degree $q^k-1$ in the functions $E_{1,v+L}$. 

For an alternative description recall that $\psi^{L\omega}_a$ can be characterised as the unique $\BF_q$-linear polynomial such that $\psi^{L\omega}_a (e_{L\omega}(z)) = e_{L\omega}(az)$. 
Plugging the expansions for $\psi^{L\omega}_a$ and $e_{L\omega}$ into this functional equation and using the fact that $e_{L,1}=1$, we deduce that for all $k\ge0$ we have
\UseTheoremCounterForNextEquation
\begin{equation}\label{CoeffsDef3}
g^L_{a,k} + \sum_{i=0}^{k-1} g^L_{a,i}\cdot e_{k-i,L}^{q^i}
\ =\ e_{k,L}\cdot a^{q^k}.
\end{equation}
By induction on~$k$ this recursion relation implies that $g^L_{a,k}$ is a universal polynomial with coefficients in $A$ in the functions $e_{j,L}$ for all $1\le j\le k$, or again in the Eisenstein series $E_{q^i-1,L}$ for all $1\le i\le k$.

\medskip
More generally, consider any non-zero ideal $N\subset A$. Then some positive power of $N$ is a principal ideal, say $N^n=(a)$ for $a\in A\setminus\{0\}$, and we choose an element $N^*\in\Cinf$ such that $(N^*)^n=a$. This element is well-defined up to multiplication by a root of unity, and for any principal ideal $(a)$ the value $(a)^*$ is equal to $a$ times a root of unity.
We also set $\deg(N) := \dim_{\BF_q}(A/N)$, so that $[N^{-1}L:L]=q^{r\deg(N)}$. In analogy with the definition (\ref{PsiLaDef}) of $\psi^{L\omega}_a$ we define 
\UseTheoremCounterForNextEquation
\begin{equation}\label{PsiLNDef}
\psi^{L\omega}_N\ :=\ N^*\cdot e_{e_{L\omega}(N^{-1}L\omega)}.
\end{equation}
Note that for any principal ideal we have $g^L_{(a),i} = g^L_{a,i}$ times a root of unity; hence everything that follows about $g^L_{N,i}$ applies equally to $g^L_{a,i}$. 

\medskip
For general $N$, by (\ref{ExpDef}) and Corollary \ref{EisensteinWeight1} we have
\UseTheoremCounterForNextEquation
\begin{equation}\label{IdealCoeffsDef1a}
\psi^{L\omega}_N(X)\ =\ 
N^*\cdot X\cdot\kern-10pt
\prod_{\textstyle{{v\,\in\,N^{-1}L\setminus L}\atop{{\rm modulo}\ L}}}
\kern-5pt \bigl(1-E_{1,v+L}(\omega)\cdot X\bigr).
\end{equation}
As in (\ref{CoeffsDef2}) we define holomorphic functions $g^L_{N,i}$ on $\Omega^r$ by expanding
\UseTheoremCounterForNextEquation
\begin{equation}\label{IdealCoeffsDef1b}
\psi^{L\omega}_N(X)\ =\ \sum_{i\ge0} g^L_{N,i}(\omega)\cdot X^{q^i},
\end{equation}
which are non-zero for $i=0$ and $i=r\deg(N)$ but zero whenever $i>r\deg(N)$. 
The formula (\ref{IdealCoeffsDef1a}) implies that each $g^L_{N,k}$ is a homogeneous symmetric polynomial of degree $q^k-1$ in the functions $E_{1,v+L}$. 

For an alternative description observe that by the definition of $\psi^{L\omega}_N$ and Proposition \ref{exp2} (a) we have 
\UseTheoremCounterForNextEquation
\begin{equation}\label{IdealCoeffsDef2}
\psi^{L\omega}_N (e_{L\omega}(z))\ =\ N^*\cdot e_{N^{-1}L\omega}(z).
\end{equation}
Plugging the respective expansions into this functional equation and using the fact that $e_{L,1}=1$, we deduce that for all $k\ge0$ we have
\UseTheoremCounterForNextEquation
\begin{equation}\label{IdealCoeffsDef3}
g^L_{N,k} + \sum_{i=0}^{k-1} g^L_{N,i}\cdot e_{k-i,L}^{q^i}
\ =\ N^*\cdot e_{k,N^{-1}L}.
\end{equation}
By induction on~$k$ this recursion relation implies that $g^L_{N,k}$ is a polynomial with coefficients in $\BF_p[N^*]$ in the functions $e_{j,L}$
and $e_{j,N^{-1}L}$ for all $1\le j\le k$, or again in the Eisenstein series $E_{q^i-1,L}$ and $E_{q^i-1,N^{-1}L}$ for all $1\le i\le k$.

\begin{Prop}\label{IdealCoeffFormIsModular}
For any non-zero ideal $N\subset A$ and any $k\ge0$ we have:
\begin{enumerate}
\item[(a)] $g^L_{N,k}|_{q^k-1}\gamma = g^{L\gamma}_{N,k}$ for all $\gamma\in\GL_r(F)$.
\item[(b)] $g^L_{N,k}$ is a modular form of weight $q^k-1$ for the group~$\Gamma_L$.
\item[(c)] The $u$-expansion of $g^L_{N,k}$ has constant term $g^{L'}_{N,k}$ with $L'$ as in (\ref{L'L1Def}).
In particular $g^L_{N,k}$ is a cusp form whenever $k>(r-1)\deg(N)$, but not for $k=(r-1)\deg(N)$.
\end{enumerate}
\end{Prop}

\begin{Proof}
By construction $g^L_{N,i}$ is a homogeneous symmetric polynomial of degree $q^i-1$ in the functions $E_{1,v+L}$. Thus the transformation formula in Proposition \ref{EisensteinBasicProps} (a) directly implies (a). Part (b) follows from Theorem \ref{EisensteinModularForm} and the formula (\ref{IdealCoeffsDef3}) by induction on~$k$.
To prove (c), write $\omega = \binom{\omega_1}{\omega'}$ as before. For any fixed $\omega'\in\Omega^{r-1}$, if $\omega_1$ goes to infinity, the defining formula (\ref{ExpDef13}) shows that $e_{L\omega}$ and $e_{N^{-1}L\omega}$ go to $e_{L'\omega'}$ and $e_{N^{-1}L'\omega'}$, respectively. The functional equation $\psi^{L\omega}_N (e_{L\omega}(z)) = N^*\cdot e_{N^{-1}L\omega}(z)$ and its counterpart for $L'\omega'$ in place of $L\omega$ thus imply that $\psi^{L\omega}_N$ goes to $\psi^{L'\omega'}_N$. Taking coefficients this shows that the $u$-expansion of each $g^L_{N,k}$ has constant term $g^{L'}_{N,k}$. Finally, that constant term is zero for $k>(r-1)\deg(N)$ and non-zero for $k=(r-1)\deg(N)$.
\end{Proof}

\section{Discriminant forms}
\label{DiscForms}

\begin{Def}\label{DiscFormDef1}
For any non-zero proper ideal $N\subsetneqq A$ we call $\Delta^L_N := g^L_{N,r\deg(N)}$ the \emph{discriminant form} associated to~$N$.
Likewise we set $\Delta^L_a := g^L_{a,r\deg(a)}$.
\end{Def}

Since $[N^{-1}L:L]$ is a power of~$q$, we have $(-1)^{[N^{-1}L:L]-1}=1$ in $\BF_q$; hence by (\ref{IdealCoeffsDef1a}) and (\ref{IdealCoeffsDef1b}) the above definition means that 
\UseTheoremCounterForNextEquation
\begin{equation}\label{DiscFormDef2}
\Delta^L_N(\omega)
\ =\ N^*\cdot\kern-10pt
\prod_{\textstyle{{v\,\in\,N^{-1}L\setminus L}\atop{{\rm modulo}\ L}}}
\kern-10pt E_{1,v+L}(\omega).
\end{equation}

\begin{Prop}\label{DiscFormProps}
\begin{enumerate}
\item[(a)] $\Delta^L_N(\omega)\not=0$ for all $\omega\in\Omega^r$. 
\item[(b)] $\Delta^L_N$ is a cusp form of weight $q^{r\deg(N)}-1$ for the group~$\Gamma_L$.
\item[(c)] $\Delta^{aL}_N = a^{1-q^{r\deg(N)}}\cdot\Delta^L_N$ for any $a\in F$.
\end{enumerate}
\end{Prop}

\begin{Proof}
(a) follows from (\ref{DiscFormDef2}) and Corollary \ref{EisensteinWeight1}, and (b) is a special case of Proposition \ref{IdealCoeffFormIsModular}. Assertion (c) results from applying Proposition \ref{IdealCoeffFormIsModular} (a) to $\gamma=a\cdot\Id_r$.
\end{Proof}

\medskip
Next recall that for any $a\in A\setminus\{0\}$ the degree $\deg(a)$ is a multiple of the degree $\deg(\infty)$  of the residue field at $\infty$ over~$\BF_q$. Therefore $q^{r\deg(a)}-1$ is a multiple of $q^{r\deg(\infty)}-1$.

\begin{Prop}\label{FundDiscForm}
There exists a non-zero cusp form $\Delta^L$ of weight $q^{r\deg(\infty)}-1$ for the group~$\Gamma_L$, such that for every $a\in A\setminus\{0\}$ we have
$$\kern30pt\Delta^L_a\ =\ (\Delta^L)^{\frac{q^{r\deg(a)}-1}{q^{r\deg(\infty)}-1}}
\cdot(\hbox{some root of unity}).$$
Moreover this $\Delta^L$ is unique up to multiplication by some root of unity.
\end{Prop}

\begin{Proof}
Since $\psi^{L}$ is a Drinfeld module, for all $a,b\in A\setminus\{0\}$ we have $\psi^L_{ab}(X) = \psi^L_a(\psi^L_b(X))$. Substituting the expansions from (\ref{CoeffsDef2}) for $\psi^L_{ab}$ and $\psi^L_a$ and $\psi^L_b$ and taking highest coefficients implies that $\Delta^L_{ab} = \Delta^L_a\cdot (\Delta^L_b)^{q^{r\deg(a)}}$. As the ring $A$ is commutative, interchanging $a$ and $b$ yields the same value; hence 
$$\Delta^L_b\cdot (\Delta^L_a)^{q^{r\deg(b)}}
\ =\ \Delta^L_a\cdot (\Delta^L_b)^{q^{r\deg(a)}}.$$
By Proposition \ref{DiscFormProps} we may divide by $\Delta^L_a\Delta^L_b$, obtaining the equality
\UseTheoremCounterForNextEquation
\begin{equation}\label{FundDiscForm1}
(\Delta^L_a)^{q^{r\deg(b)}-1}\ =\ (\Delta^L_b)^{q^{r\deg(a)}-1}.
\end{equation}

To exploit this fact, recall that by the Riemann-Roch theorem, every sufficiently large multiple of $\deg(\infty)$ arises as $\deg(a)$ for some element $a\in A\setminus\{0\}$. In particular we can find non-constant elements $b$, $c\in A$ such that $\deg(b)=\deg(c)+\deg(\infty)$. Then by Proposition \ref{DiscFormProps} the quotient
\UseTheoremCounterForNextEquation
\begin{equation}\label{FundDiscForm2}
\Delta^L\ := \Delta^L_b \big/(\Delta^L_c)^{q^{r\deg(\infty)}}
\end{equation}
is a well-defined holomorphic function on~$\Omega^r$. The fact that $\Delta^L_b$ and $\Delta^L_c$ are modular forms of respective weights $q^{r\deg(b)}-1$ and $q^{r\deg(c)}-1$ for $\Gamma_L$ implies that $\Delta^L$ is a weak modular form of weight 
$$\bigl(q^{r\deg(b)}-1\bigr) - 
\bigl(q^{r\deg(c)}-1\bigr)\cdot q^{r\deg(\infty)}
\ =\ q^{r\deg(\infty)}-1$$
for~$\Gamma_L$. Also, by direct calculation the formula (\ref{FundDiscForm1}) in the case $a=c$ implies that 
$$(\Delta^L)^{q^{r\deg(b)}-1}\ =\ (\Delta^L_b)^{q^{r\deg(\infty)}-1}.$$
Combining this with the formula (\ref{FundDiscForm1}) for arbitrary $a$
we deduce that
$$(\Delta^L_a)^{(q^{r\deg(\infty)}-1)(q^{r\deg(b)}-1)}
\ =\ (\Delta^L)^{(q^{r\deg(a)}-1)(q^{r\deg(b)}-1)}.$$
Thus $\Delta^L_a / (\Delta^L)^{\frac{q^{r\deg(a)}-1}{q^{r\deg(\infty)}-1}}$ is a holomorphic function on $\Omega^r$ whose $(q^{r\deg(\infty)}-1)(q^{r\deg(b)}-1)$-th power is identically~$1$. As the rigid analytic space $\Omega^r$ is connected, this function is therefore constant and a root of unity. The last formula also shows that a positive power of $\Delta^L$ is holomorphic at every boundary component; hence the same holds for~$\Delta^L$. Thus $\Delta^L$ has all the desired properties. Finally, the uniqueness is clear from the stated condition.
\end{Proof}

\begin{Prop}\label{FundDiscFormN}
For every non-zero proper ideal $N\subsetneqq A$ we have
$$\kern30pt\Delta^{N^{-1}L}\cdot(\Delta^L_N)^{q^{r\deg(\infty)}-1}
\ =\ (\Delta^L)^{q^{r\deg(N)}-1}\cdot(\hbox{some constant}).$$
\end{Prop}

\begin{Proof}
The formulas (\ref{PsiLNDef}) and (\ref{IdealCoeffsDef2}) imply that $\psi^{L\omega}_N = N^*\cdot h^L_N$, where $h^L_N$ is an isogeny of Drinfeld modules $(\BGacomma{\Omega^r},\psi^{L}) \to (\BGacomma{\Omega^r},\psi^{N^{-1}L})$. 
For any $a\in A$ we then have $\psi^{N^{-1}L}_a\circ h^L_N = h^L_N\circ\psi^L_a$. Taking highest coefficients implies that 
$$\Delta^{N^{-1}L}_a \cdot (\Delta^L_N)^{q^{r\deg(a)}} 
\ =\ \Delta^L_N \cdot (\Delta^L_a)^{q^{r\deg(N)}} \cdot (\hbox{some constant}).$$
Dividing by $\Delta^L_N$ and substituting the formulas for $\Delta^{N^{-1}L}_a$ and $\Delta^L_a$ from Proposition \ref{FundDiscForm} we obtain
$$(\Delta^{N^{-1}L})^{\frac{q^{r\deg(a)}-1}{q^{r\deg(\infty)}-1}} \cdot (\Delta^L_N)^{q^{r\deg(a)}-1} 
\ =\ (\Delta^L)^{q^{r\deg(N)}\cdot\frac{q^{r\deg(a)}-1}{q^{r\deg(\infty)}-1}}
\cdot (\hbox{some constant}).$$
Varying $a$ or extracting roots as in Proposition \ref{FundDiscForm} yields the desired formula.
\end{Proof}


\begin{Rem}\label{FundDiscRem}
\rm If the class group $\Cl(A)$ of $A$ is trivial, the above relations show that $\Delta^L$ is the unique fundamental discriminant form for~$\Gamma_L$.

In general, for any non-zero proper ideal $M\subsetneqq A$ we have $\Gamma_{M^{-1}L}=\Gamma_L$. The discriminant forms $\Delta^{M^{-1}L}_a$ and $\Delta^{M^{-1}L}$ and $\Delta^{M^{-1}L}_N$ are therefore cusp forms for the same group~$\Gamma_L$. 
Let $\CH$ denote the multiplicative group of nowhere vanishing holomorphic functions on $\Omega^r$ up to constants that is generated by all of them. Then the formulas in Propositions \ref{DiscFormProps} (c) and \ref{FundDiscForm} and \ref{FundDiscFormN} imply that as $N$ runs through a set of representatives of the ideal class group $\Cl(A)$, the functions $\Delta^{N^{-1}L}$ generate a subgroup of finite index, say~$\CH'$.

On the other hand each discriminant form corresponds to a section of a certain invertible sheaf on the Satake compactification of $\Gamma_L\backslash\Omega$. As such, its divisor is a formal $\BZ$-linear combination of the irreducible components of codimension $1$ of the boundary. These irreducible components are  in bijection with $\Cl(A)$, so the group $\CD$ of divisors supported on the boundary is a free abelian group of rank $\Cl(A)$. Taking divisors maps the above group $\CH$ injectively into~$\CD$.

One can expect that the image of $\CH$ has finite index in~$\CD$. In fact, precisely such a statement is proved for an arbitrary congruence subgroup in the case $r=2$ by Gekeler \cite[VII Thm.\;5.11]{GekelerDMC} and \cite[Thm.\;4.1]{Gekeler1997}, 
and by Kapranov \cite[top of page 546]{Kapranov} for arbitrary $r$ in the case $A=\BF_q[t]$.

Note that, since $\CH'$ is generated by $|\Cl(A)|$ elements and has finite index in~$\CH$, the expectation is equivalent to saying that $\CH$ is a free abelian group of rank $|\Cl(A)|$. This in turn means that, up to taking roots, the formulas in Propositions \ref{DiscFormProps} (c) and \ref{FundDiscForm} and \ref{FundDiscFormN} generate all multiplicative relations up to constant factors between the discriminant forms.
\end{Rem}

\begin{Ex}\label{FundDiscEx}
\rm Suppose that $\Spec A$ is a rational curve and $\infty$ is a point of degree $2$ over~$\BF_q$. Let $P\subset A$ be the prime ideal associated to a point of degree $1$ over~$\BF_q$. Then the ideal class group of $A$ has order $2$ and is generated by the class of~$P$. Write $P^2=(a)$ for an element $a\in A$ of degree~$2$. Then by Proposition \ref{FundDiscForm} we have $\Delta^L_a\sim\Delta^L$, where ``$\sim$'' denotes equality up to a constant. 
Also, in the notation of the proof of Proposition \ref{FundDiscFormN} we have $a\cdot h^{P^{-1}L}_P\circ h^L_P = \psi^L_a$. Taking highest coefficients implies that $\Delta^{P^{-1}L}_P\cdot (\Delta^L_P)^{q^r} \sim \Delta^L_a \sim \Delta^L$. Together with the same relation for $P^{-1}L$ in place of~$L$ and with the fact that $\Delta^{P^{-2}L}_P = \Delta^{a^{-1}L}_P \sim \Delta^L_P$ by Proposition \ref{DiscFormProps} (c), we conclude that
$$\begin{array}{rl}
\Delta^{P^{-1}L}_P\cdot (\Delta^L_P)^{q^r} &\!\!\sim\ \Delta^L \quad\hbox{and}\\[3pt]
\Delta^L_P\cdot (\Delta^{P^{-1}L}_P)^{q^r} &\!\!\sim\ \Delta^{P^{-1}L}.
\end{array}$$
In this case we can therefore view $\Delta^L_P$ and $\Delta^{P^{-1}L}_P$ as the two fundamental discriminant forms for~$\Gamma_P$, and by Remark \ref{FundDiscRem} they should be multiplicatively independent.
\end{Ex}


\begin{Rem}\label{DeltaProductExpansion}
\rm In the case $A=\BF_q[t]$ one can take $\Delta^L = \Delta^L_t$ in Proposition \ref{FundDiscForm}. In \cite{basson2016} this function is shown to satisfy a product formula which generalises the Jacobi product formula in the rank 2 case of Gekeler \cite{GekelerProd}. Another product formula, involving $r-1$ separate
parameters with constant coefficients, rather than $u$-expansions treated in the present paper, was obtained by Hamahata \cite{Hamahata}.
\end{Rem}

\begin{Rem}\label{NonZeroCuspForm}
\rm For any $v\in F^r\setminus L$, the Eisenstein series $E_{1,v+L}$ is a non-zero modular form of weight $1$ for the group $\Gamma_{v+L}$ by Corollary \ref{EisensteinWeight1} and Theorem \ref{EisensteinModularForm}. Using Proposition \ref{FundDiscForm} it follows that for any integer $k\ge0$, the product $\Delta^L\cdot E_{1,v+L}^k$ is a non-zero cusp form of weight $q^{r\deg(\infty)}-1+k$ for $\Gamma_{v+L}$. In this way we can explicitly produce non-zero cusp forms for $\Gamma_{v+L}$ of any sufficiently large weight, giving more substance to the abstract result of Proposition \ref{CuspEx}.
\end{Rem}


To finish this section we construct Drinfeld modular forms of non-zero type by extracting roots from discriminant forms.
This rests on the observation that for every $\alpha\in\BF_q^\times$, applying Proposition \ref{EisensteinBasicProps} (a) to $\gamma=\alpha\cdot\Id_r$ implies that 
\UseTheoremCounterForNextEquation
\begin{equation}\label{AlphaTrafo}
E_{1,\alpha v+L}\ =\ E_{1,v+L}\,|_1\,\alpha\cdot\Id_r\ 
\stackrel{\eqref{IdrAction}}{=}\ \alpha^{-1}\cdot E_{1,v+L}.
\end{equation}
Plugging this into (\ref{DiscFormDef2}), we can write each discriminant form as a $(q-1)$-st power of another holomorphic function on~$\Omega^r$. 

Specifically, choose a set of representatives $\CR^L_N$ of $N^{-1}L\setminus L$ modulo addition by $L$ and multiplication by~$\BF_q^\times$. Choose an element $\lambda_N\in\Cinf$ satisfying $\lambda_N^{q-1} = -N^*$.
Consider the function
\UseTheoremCounterForNextEquation
\begin{equation}\label{RootDiscFormDef}
\delta^L_N(\omega)
\ :=\ \lambda_N \cdot \prod_{v\in\CR^L_N} E_{1,v+L}(\omega).
\end{equation}

\begin{Prop}\label{RootDiscFormProps}
\begin{enumerate}
\item[(a)] We have $(\delta^L_N)^{q-1} = \Delta^L_N$. In particular, another choice of representatives or of $\lambda_N$ changes $\delta^L_N$ only by a factor in~$\BF_q^\times$.
\item[(b)] The function $\delta^L_N$ is a cusp form of weight $\frac{q^{r\deg(N)}-1}{q-1}$ and type $\deg(N)$ for the group~$\Gamma_L$.
\end{enumerate}
\end{Prop}

\begin{Proof}
Abbreviate $k:= |\CR^L_N| = \frac{q^{r\deg(N)}-1}{q-1}$ and note that $(\prod_{\alpha\in\BF_q^\times}\alpha)^k=(-1)^k=-1$. Using this, the definitions of $\delta^N_L$ and $\Delta^N_L$ and (\ref{AlphaTrafo}) imply that
$$(\delta^L_N)^{q-1}
\ =\ - N^*\cdot \kern-4pt\prod_{v\in\CR^L_N} E_{1,v+L}^{q-1}
\ =\ N^*\cdot \kern-4pt\prod_{v\in\CR^L_N} \prod_{\alpha\in\BF_q^\times}\alpha^{-1}E_{1,v+L}
\ =\ N^*\cdot \kern-4pt\prod_{v\in\CR^L_N} \prod_{\alpha\in\BF_q^\times}E_{1,\alpha v+L}
\ =\ \Delta^L_N(\omega),$$
proving (a). The proof of (b) rests on properties of the Moore determinant, assembled in \cite[Chapter 1.3]{GossBS}: For any elements  $x_1,\ldots,x_n$ of an $\BF_q$-algebra the Moore determinant is defined as
\UseTheoremCounterForNextEquation
\begin{equation}\label{Moore1}
M(x_1,x_2,\ldots,x_n)\ :=\ \left\lvert
\begin{array}{ccc} 
x_1 & \cdots & x_n \\ 
x_1^q & & x_n^q \\ 
\vdots & & \vdots \\ 
x_1^{q^{n-1}}\kern-4pt & \cdots & x_n^{q^{n-1}}\kern-4pt
\end{array}\right\rvert.
\end{equation}
Its most important property is \cite[Cor.\;1.3.7]{GossBS}
\UseTheoremCounterForNextEquation
\begin{equation}\label{Moore2}
M(x_1,x_2,\ldots,x_n)\ =\ \prod_{(\alpha_1,\ldots,\alpha_n)}
\Bigl(\sum_{i=1}^n \alpha_i x_i \Bigr),
\end{equation}
where the product extends over all tuples in $\BF_q^n\setminus\{(0,\ldots,0)\}$ whose first non-zero entry is~$1$.
Also, for any matrix $B=(\beta_{ij})_{i,j=1,\ldots,n}$ with coefficients in $\BF_q$ we have $\beta_{ij}^q=\beta_{ij}$; hence the multiplicativity of the determinant implies that
\UseTheoremCounterForNextEquation
\begin{equation}\label{Moore3}
 M\Bigl(\sum_{i=1}^n\beta_{i1}x_i, \ldots, \sum_{i=1}^n\beta_{in}x_i\Bigr) 
\ =\ \det(B)\cdot M(x_1,x_2,\ldots,x_n).
\end{equation}

To apply this, choose elements $v_1,\ldots,v_n\in N^{-1}L\setminus L$ whose residue classes form a basis of the $\BF_q$-vector space $N^{-1}L/L$. Then the set $\CR^L_N$ of all elements of the form $\sum_{i=1}^n \alpha_i v_i$, for tuples $(\alpha_1,\ldots,\alpha_n)\in\BF_q^n\setminus\{(0,\ldots,0)\}$ whose first non-zero entry is~$1$, is a set of representatives of $N^{-1}L\setminus L$ modulo addition by $L$ and multiplication by~$\BF_q^\times$. The formula (\ref{Moore2}) and the additivity of the exponential function then imply that
$$M\bigl(e_{L\omega}(v_1\omega),\ldots,e_{L\omega}(v_n\omega)\bigr)
\ = \kern-5pt \prod_{(\alpha_1,\ldots,\alpha_n)}
\Bigl(\sum_{i=1}^n \alpha_i e_{L\omega}(v_i\omega) \Bigr)
\ =\ \prod_{v\in\CR^L_N} e_{L\omega}(v\omega).$$
Take an arbitrary element $\gamma\in\Gamma_L$. Then the same calculation with the basis $v_1\gamma,\ldots,v_n\gamma$ yields
$$M\bigl(e_{L\omega}(v_1\gamma\omega),\ldots,e_{L\omega}(v_n\gamma\omega)\bigr)
\ =\ \prod_{v\in\CR^L_N} e_{L\omega}(v\gamma\omega).$$
For each $j$ choose $\beta_{ij}\in\BF_q$ such that $v_j\gamma \equiv \sum_{i=1}^n\beta_{ij}v_i$ modulo~$L$. Then by the $\BF_q$-linearity of the exponential function we have $e_{L\omega}(v_j\gamma\omega) = \sum_{i=1}^n\beta_{ij}e_{L\omega}(v_i\gamma\omega)$; hence with $B:=(\beta_{ij})_{i,j=1,\ldots,n}$ the formula (\ref{Moore3}) implies that 
$$M\bigl(e_{L\omega}(v_1\gamma\omega),\ldots,e_{L\omega}(v_n\gamma\omega)\bigr)
\ =\ \det(B)\cdot M\bigl(e_{L\omega}(v_1\omega),\ldots,e_{L\omega}(v_n\omega)\bigr).$$
Combining these computations we deduce that 
$$\prod_{v\in\CR^L_N} e_{L\omega}(v\gamma\omega)
\ =\ \det(B)\cdot\prod_{v\in\CR^L_N} e_{L\omega}(v\omega).$$
Using Proposition \ref{EisensteinBasicProps} (a) and Corollary \ref{EisensteinWeight1} we find that 
\begin{eqnarray*}
(\delta^L_N|_k\gamma)(\omega)
\!\!&=&\!\! \lambda_N \cdot \prod_{v\in\CR^L_N} (E_{1,v+L}|_1\gamma)(\omega)
\ =\ \lambda_N \cdot \prod_{v\in\CR^L_N} E_{1,v\gamma+L}(\omega) \\
\!\!&=&\!\! \lambda_N \cdot \det(B)^{-1} \cdot \prod_{v\in\CR^L_N} E_{1,v+L}(\omega)
\ =\ \det(B)^{-1} \cdot \delta^L_N(\omega).
\end{eqnarray*}
To determine $\det(B)$ note that since $L$ is a projective module of rank $r$ over~$A$, the module $N^{-1}L/L$ is a free module of rank $r$ over $A/N$. Without loss of generality we may therefore assume that the $\BF_q$-basis $v_1,\ldots,v_n$ is formed by multiplying an $A/N$-basis of $N^{-1}L/L$ with an $\BF_q$-basis of $A/N$. For a suitable order of this basis, the matrix $A$ is then simply a block diagonal matrix with $m := \dim_{\BF_q}(A/N) = \deg(N)$ copies of $\gamma$ on the diagonal. Therefore $\det(B)=\det(\gamma)^m$. In view of (\ref{DefOfActionOnF}) the above calculation thus implies that 
$$\delta^L_N|_{k,m}\gamma
\ =\ \det(\gamma)^{m}\cdot\delta^L_N|_k\gamma
\ =\ \delta^L_N.$$
In other words $\delta^L_N$ is a weak modular form of weight $k$ and type $m$ for the group~$\Gamma_L$. But by Theorem \ref{EisensteinModularForm} and construction it is already a modular form for the congruence subgroup $\Gamma_L(N)$. It is therefore a modular form for~$\Gamma_L$. Finally, since $\Delta^L_N$ is a cusp form, assertion (a) implies that $\delta^L_N$ is a cusp form as well. This finishes the proof of (b).
\end{Proof}


\begin{Rem}\rm
In the case $A=\BF_q[t]$ and $L=A^r$, the cusp form $\delta_t^L$ was first constructed by Gekeler in the 1980's, and is called $h(\omega)$ in the literature. 
The $r=2$ case appears in \cite{GekelerCoeff} while the $r>2$ case was unpublished until \cite{GekelerHigherRankI}. In the meantime, it made an appearance as a weak modular form in \cite{GekelerQuasi} and was shown to be holomorphic at infinity by Perkins \cite{Perkins}. 
In \cite[Thm.\,5.3]{BassonBreuer2017} it is shown to satisfy a product formula derived from the product formula of $\Delta_t^L$.
\end{Rem}


\section{The special case $A=\BF_q[t]$}
\label{SecFqt}

Throughout this section we set $A:=\BF_q[t]$ and $L:=A^r$. Then $\Gamma_L = \GL_r(A)$, and $\Gamma(t):=\Gamma_L((t))$ is the subgroup of matrices in $\GL_r(A)$ which are congruent to the identity matrix modulo $(t)$. 
Recall from (\ref{ModRing}) that the graded ring of modular forms of all weights for an arithmetic group $\Gamma$ is defined as
$$\CM_{*}(\Gamma)\ :=\ \bigoplus_{k\ge0}\CM_{k}(\Gamma).$$
For $\Gamma=\Gamma(t)$ this ring can be described very explicitly, and for a subgroup containing $\Gamma(t)$ a description can be deduced by taking invariants.
In the case $r=2$ the ring was determined by Cornelissen \cite{cornelissen_drinfeld_1997-1} for $\Gamma(t)$, by Goss \cite{goss_modular_1980} for $\GL_2(A)$, and by Gekeler \cite{GekelerCoeff} for $\SL_2(A)$.


\begin{Thm}\label{ModRing1}
The ring $\CM_{*}(\Gamma(t))$ is generated over $\Cinf$ by the Eisenstein series $E_{1,v+L}$ of weight $1$ for all $v\in t^{-1}L\setminus L$, and all polynomial equations between them are induced by the relations
$$\begin{array}{rll}
E_{1,\alpha v+L} &\!\!=\ \alpha^{-1}\cdot E_{1,v+L}
& \mbox{for all $v \in t^{-1}L\setminus L$ and $\alpha \in \BF_q^{\times}$, and} \\[5pt]
E_{1,v+L}\cdot E_{1,v'+L} &\!\!=\ E_{1,v+v'+L}\cdot (E_{1,v+L}+E_{1,v'+L})\ 
& \mbox{for all $v$, $v' \in t^{-1}L\setminus L$ with $v+v' \not\in L$.}
\end{array}$$
\end{Thm}

\begin{Proof}
Let $K(t)<\GL_r(\hatA)$ denote the subgroup of matrices that are congruent to the identity matrix modulo $(t)$. By construction it is open compact and fine in the sense of \cite[\DefinitionCite1.4]{Pink}. Let $\smash{M^r_{A,K(t)}}$ be the associated fine moduli space of Drinfeld $A$-modules of rank $r$ with a full level $(t)$ structure. Then $\GL_r(\BAfF) = \GL_r(F)\cdot K(t)$, and so (\ref{AnalUnivPig1b2}) with $g=1$ provides an isomorphism $\pi_1: \Gamma(t)\backslash\Omega^r \stackrel{\sim}{\longto} \smash{M^r_{A,K(t)}}(\Cinf)$. The Satake compactification $\smash{\OM^r_{A,K(t)}}$ was described explicitly in \cite{PinkSchieder} and \cite{Pink}, as follows.

Abbreviate $\bar V := t^{-1}L/L$, and let $A_{\bar V}$ denote the graded polynomial ring over $\BF_q$ in independent variables $Y_{\bar v}$ of degree $1$ for all $\bar v\in \bar V\setminus\{0\}$. Let $\Fa_{\bar V} \subset A_{\bar V}$ be the homogeneous ideal that is generated by the elements of the form
$$\begin{array}{cl}
Y_{\alpha\bar v} - \alpha^{-1} Y_{\bar v}
& \mbox{for all $\bar v \in \bar V\setminus\{0\}$ and $\alpha \in \BF_q^{\times}$, and} \\[5pt]
Y_{\bar v} Y_{\bar v'} - Y_{\bar v+\bar v'} \cdot (Y_{\bar v}+Y_{\bar v'})\ 
& \mbox{for all $\bar v$, $\bar v' \in \bar V\setminus\{0\}$ with $\bar v + \bar v'\not=0$.}
\end{array}$$
Let $R_{\bar V} := A_{\bar V}/\Fa_{\bar V}$ denote the graded factor ring. Then by \cite[Thm.\;7.4]{Pink} there is a natural isomorphism 
\UseTheoremCounterForNextEquation
\begin{equation}\label{ModRing1a}
\OM^r_{A,K(t)}\ \cong\ \Proj (R_{\bar V}\otimes_{\BF_q} F),
\end{equation}
which also identifies the invertible sheaf $\CL$ from Section \ref{Sec:AAMF} with the ample sheaf $\CO(1)$ on $\Proj (R_{\bar V}\otimes_{\BF_q} F)$. Combined with Theorem \ref{AAMFThm1} 
we thus obtain an isomorphism of graded $\Cinf$-algebras 
\UseTheoremCounterForNextEquation
\begin{equation}\label{ModRing1b}
\CM_{*}(\Gamma(t))\ \cong\ R_{\bar V}\otimes_{\BF_q} \Cinf.
\end{equation}
By the proof of \cite[Thm.\;7.4]{Pink}, the isomorphism (\ref{ModRing1a}) also realises the universal generalised Drinfeld $A$-module over \smash{$\OM^r_{A,K(t)}$} as the pair $(\bar{E},\bar{\phi})$ consisting of the line bundle whose sheaf of sections is the invertible sheaf dual to $\CO(1)$ and the generalised Drinfeld $A$-module with
$$\bar\phi_t(X)\ =\ t\cdot X\cdot \prod_{\bar v\in\bar V\setminus\{0\}} (1-\bar Y_{\bar v}\cdot X),$$
where $\bar Y_{\bar v}\in R_V$ denotes the residue class of~$Y_{\bar v}$. 
On the other hand from (\ref{StoneTrek}) we have a natural isomorphism
$$\pi_g^*(\bar E,\bar\phi) \ \cong\ (\BGacomma{\Omega^r},\psi^{L}),$$
and by equation (\ref{CoeffsDef1}) we have
$$\psi^{L\omega}_t(X)
\ =\ t\cdot X\cdot\kern-10pt
\prod_{\textstyle{{v\,\in\,t^{-1}L\setminus L}\atop{{\rm modulo}\ L}}}
\kern-5pt \bigl(1-E_{1,v+L}(\omega)\cdot X\bigr).$$
Furthermore, the respective level structures send a non-zero residue class $\bar v=v+L$ to the element $\bar Y_{\bar v}^{-1}$ in one case and to the function $E_{1,v+L}(\omega)^{-1} = e_{L\omega}(v\omega)$ in the other. Under the isomorphism (\ref{ModRing1b}) the element $\bar Y_{\bar v}$ therefore corresponds precisely to the Eisenstein series $E_{1,v+L}$. By the construction of $R_{\bar V}$ these Eisenstein series therefore generate $\CM_{*}(\Gamma(t))$ and satisfy precisely the stated algebraic relations.
\end{Proof}

\begin{Cor}\label{ModRing1Quot}
The quotient field of $\CM_{*}(\Gamma(t))$ is a rational function field over $\Cinf$ that is generated by the algebraically independent elements $E_{1,v_i+L}$ as $v_i+L$ runs through any $\BF_q$-basis of $t^{-1}L/L$.
\end{Cor}

\begin{Proof}
By \cite{PinkSchieder} the ring $R_V$ is an integral domain and its quotient field is a rational function field over $\BF_q$ that is generated by the algebraically independent elements $\bar Y_{\bar v_i}$ for any basis $\bar v_1,\ldots,\bar v_r$ of~$\bar V$. The corollary thus follows from the isomorphism (\ref{ModRing1b}).
\end{Proof}


\begin{Thm}\label{ModRing2}
\begin{enumerate}
\item[(a)] The ring $\CM_{*}(\GL_r(A))$ is generated over $\Cinf$ by the coefficient forms $g^L_{t,i}$ of weight $q^i-1$ for all $1\le i\le r$, which are algebraically independent over~$\Cinf$. The same statement holds with the coefficient forms $e_{i,L}$ or the Eisenstein series $E_{q^i-1,L}$ in place of~$g^L_{t,i}$.
\item[(b)] The ring $\CM_{*}(\SL_r(A))$ is generated over $\Cinf$ by the coefficient forms $g^L_{t,i}$ of weight $q^i-1$ for all $1\le i\le r-1$ and the determinant form $\delta^L_t$ of weight $\frac{q^r-1}{q-1}$, which are algebraically independent over~$\Cinf$. The same statement holds with the coefficient forms $e_{i,L}$ or the Eisenstein series $E_{q^i-1,L}$ in place of~$g^L_{t,i}$.
\item[(c)] Let $\Gamma_1(t)$ denote the subgroup of matrices in $\GL_r(A)$ which are congruent modulo $(t)$ to an upper triangular matrix with diagonal entries~$1$. The ring $\CM_{*}(\Gamma_1(t))$ is generated over $\Cinf$ by the modular forms
$$\sum_{\alpha_{i+1},\ldots,\alpha_r\in\BF_q} \kern-5pt
E_{1,t^{-1}(0,\ldots,0,1,\alpha_{i+1},\ldots,\alpha_r)+L}$$
of weight $1$ for all $1\le i\le r$, which are algebraically independent over~$\Cinf$.
\end{enumerate}
\end{Thm}

\begin{Proof}
For any subgroup $\Gamma<\GL_r(A)$ containing $\Gamma(t)$, the formula (\ref{ModInvs}) shows that $\CM_{*}(\Gamma)$ is the subring of $\Gamma$-invariants in $\CM_{*}(\Gamma(t))$ for the natural action by $f\mapsto f|_k\gamma$ on each $\CM_k(\Gamma(t))$. By Proposition \ref{EisensteinBasicProps} (a) the action is given on the generators of $\CM_{*}(\Gamma(t))$ by $E_{1,v+L}|_1\gamma = E_{1,v\gamma+L}$. This action factors through the factor group $\Gamma/\Gamma(t)$, which is $\GL_r(\BF_q)$ in the case (a), respectively $\SL_r(\BF_q)$ in the case (b), respectively the subgroup of upper triangular matrices with diagonal entries~$1$ in the case (c). Using a theorem of Dickson, the respective ring of invariants was shown in \cite[Theorem 3.1]{PinkSchieder} to have the set of generators that is first named in each case.
The recursion relations (\ref{CoeffsDef3}) and (\ref{ExpCoeff1a}) imply that by induction on~$i$, each generator $g^L_{a,i}$ can be replaced by $e_{i,L}$ or again by $E_{q^i-1,L}$.

Since we are taking invariants under a finite group, the ring $\CM_{*}(\Gamma(t))$ is an integral extension of $\CM_{*}(\Gamma(t))^\Gamma$. The respective quotient fields therefore have the same transcendence degree over~$\Cinf$. For the former this transcendence degree is $r$ by Corollary \ref{ModRing1Quot}. In each case the $r$ given generators of the subring $\CM_{*}(\Gamma(t))^\Gamma$ must therefore be algebraically independent over~$\Cinf$.
\end{Proof}

\begin{Thm}\label{ModRing3}
For any integer $k$ we have 
$$\CM_k(\SL_r(A))\ = \bigoplus_{0\le m< q-1}\kern-5pt \CM_{k,m}(\GL_r(A)).$$
In addition, for any integer $0\le m<q-1$ we have
$$\CM_{k,m}(\GL_r(A))\ =\ (\delta^L_t)^m \cdot \CM_{k-m\frac{q^r-1}{q-1}}(\GL_r(A)).$$
In particular, every modular form for $\GL_r(A)$ of type $\not\equiv0$ modulo $(q-1)$ is a cusp form.
\end{Thm}

\begin{Proof}
The determinant induces an isomorphism $\GL_r(\BF_q)/\SL_r(\BF_q) \stackrel{\sim}{\longto} \BF_q^\times$; hence the action $f\mapsto f|_k\gamma$ of $\GL_r(\BF_q)$ on $\CM_k(\SL_r(A))$ factors through an action of~$\BF_q^\times$. As any linear action of $\BF_q^\times$ on an $\BF_q$-vector space is diagonalisable, it follows that $\CM_k(\SL_r(A))$ is a direct sum of eigenspaces. By Definition \ref{Def:WeakModForm} and (\ref{DefOfActionOnF}) these eigenspaces are just the spaces $\CM_{k,m}(\GL_r(A))$, proving the first equality.

The descriptions from Theorem \ref{ModRing2} (a) and (b) imply that $\CM_{*}(\SL_r(A))$ is a free module with basis 
$1,\delta^L_t,\ldots,(\delta^L_t)^{q-2}$
over the subring $\CM_{*}(\GL_r(A))$. Since $(\delta^L_t)^m$ is a modular form of weight $m\frac{q^r-1}{q-1}$, this results in the second assertion. 
The last one now follows from the fact that $\delta^L_t$ is a cusp form.
\end{Proof}

\begin{Rem}\label{ModRing3Explan}
\rm The last statement of Theorem \ref{ModRing3} was already established independently in Corollary \ref{Cor:GLrACase} (b) using the $u$-expansion. Combined with Proposition \ref{DiscVanishingOrder} below and the fact that $\delta^L_t$ is a modular form of weight $\frac{q^r-1}{q-1}$ and type $1$ it directly implies the second statement of Theorem \ref{ModRing3} by induction on~$m$.
\end{Rem}

\begin{Prop}\label{DiscVanishingOrder}
The Satake compactification $\OM^r_{A,\GL_r(A)}$ has only one boundary component of codimension~$1$, and the cusp form $\delta^L_t$ has vanishing order $1$ there.
\end{Prop}

\begin{Proof}
The first statement can be deduced from the fact from Proposition \ref{Prop:cosets} (a) that $\GL_r(\BAfF) = \GL_r(A)\cdot P(F)$ with the parabolic subgroup $P<\GL_r$ from (\ref{DefPL}).

For the second statement note first that under the isomorphism $\iota$ of (\ref{eq:iota}) the subgroups ${\Gamma(t)\cap U(F)} \allowbreak < \GL_r(A)\cap U(F)$ correspond to the subgroups $(At)^{r-1} \subset A^{r-1}$ of $F^{r-1}$, which have index $q^{r-1}$ in each other. Now consider any element $v\in t^{-1}L\setminus L$. By the proof of Proposition \ref{EisensteinVanishingOrderEx} the subgroup ${\Gamma_{v+L}\cap U(F)}$ corresponds to the subgroup $(At)^{r-1}$ if $v\not\in L+(\{0\}\times F^{r-1})$. By Proposition \ref{EisensteinVanishingOrderEx} we thus have 
$$\ord_{\Gamma(t)\cap U(F)}(E_{1,v+L})
\ =\ \ord_{\Gamma_{v+L}\cap U(F)}(E_{1,v+L})
\ =\ \biggl\{\begin{array}{ll}
0 & \hbox{if $v\in L+(\{0\}\times F^{r-1})$,}\\
1 & \hbox{otherwise.}
\end{array}$$
Taking the product over a set of representatives as in (\ref{RootDiscFormDef}), where the second case occurs $\frac{q^r-q^{r-1}}{q-1}=q^{r-1}$ times, we deduce that 
$$\ord_{\Gamma(t)\cap U(F)}(\delta^L_t) \ =\ q^{r-1}.$$
Since $[\GL_r(A)\cap U(F):\Gamma(t)\cap U(F)]=q^{r-1}$, it follows that $\ord_{\GL_r(A)\cap U(F)}(\delta^L_t) = 1$, as desired.
\end{Proof}


\begin{Cor}\label{ModRing2cusp1}
The cusp forms of all weights and type $0$ for $\GL_r(A)$ form the principal ideal of $\CM_{*}(\GL_r(A))$ that is generated by $\Delta_t^L$. In other words, for every integer $k$ we have 
$$\CS_{k}(\GL_r(A))\ =\ \Delta_t^L\cdot \CM_{k-q^r+1}(\GL_r(A)).$$
\end{Cor}

\begin{Proof}
The cusp form $\delta^L_t$ is non-zero everywhere by Propositions \ref{DiscFormProps} (a) and \ref{RootDiscFormProps} (a). Thus for every cusp form $f\in \CS_{k,0}(\Gamma)$, the quotient $f/\delta^L_t$ is again a weak modular form, and by Proposition \ref{DiscVanishingOrder} it is holomorphic at infinity; hence $f/\delta^L_t \in \CM_{k-\frac{q^r-1}{q-1},-1}(\GL_r(A))$. By Theorem \ref{ModRing3} with $m=q-2$ this in turn implies that $f \in (\delta^L_t)^{q-1}\CM_{k-q^r+1,0}(\GL_r(A))$, as desired.
\end{Proof}

\begin{Cor}\label{ModRing2cusp2}
The space of cusp forms $\CS_{k}(\GL_r(A))$ is zero for $k<q^r-1$ and one-dimensional with basis $\Delta_t^L$ for $k=q^r-1$. In particular $\Delta_t^L$ is an eigenform for the Hecke operator associated to any double coset $\GL_r(A)\delta\GL_r(A)\subset\GL_r(F)$.
\end{Cor}

\begin{Proof}
By Theorem \ref{ModRing2} (a) we have $\CM_{k}(\GL_r(A))=0$ for $k<0$ and $=\Cinf$ for $k=1$. By Corollary \ref{ModRing2cusp1} this implies the first statement, which in turn implies the second.
\end{Proof}


\begin{Thm}\label{DimensionFormulas}
We have the following dimension formulas for all $k\ge0$ and $m$:
\begin{enumerate}
\item[(a)] $\displaystyle \dim_{\Cinf}\MM_k(\Gamma(t)) =
\kern-10pt \sum_{i_1, \ldots, i_{r-1} \in \{0,1\}} \kern-10pt q^{\sum_{\nu} \nu \cdot i_{\nu}} \cdot \binom{k}{\sum_{\nu} i_{\nu}}$.
\item[(b)] Denote by $P_S(k)$ the number of partitions of $k$ with parts in $S=\{q-1,q^2-1,\ldots,q^r-1\}$. Then

$\displaystyle \dim_{\Cinf}\MM_k(\GL_r(A)) = P_S(k)
 = \left\{\begin{array}{ll} 0 & \text{if $(q-1)\nmid k$,}\\[3pt]
\displaystyle \frac{1}{\prod_{i=2}^r(q^i-1)} \cdot \frac{k^{r-1}}{(r-1)!} 
+ O(k^{r-2}) & \text{if $(q-1)|k$.} \end{array}\right.$
\item[(c)] $\dim_{\Cinf}\MM_{k,m}(\GL_r(A)) = 
\left\{\begin{array}{ll} 
P_S\big(k-m\frac{q^r-1}{q-1}\big) & \text{if $k\ge m\frac{q^r-1}{q-1}$,}\\[3pt]
0 & \text{otherwise.} \end{array}\right.$
\item[(d)] $\dim_{\Cinf}\MM_k(\Gamma_1(t)) = \binom{k-1}{r-1}$.
\end{enumerate}
\end{Thm}

\begin{Proof}
Assertion (a) follows from Theorem \ref{ModRing1} together with \cite[Thm.\;1.10]{PinkSchieder}. 
The first equality in (b) results from \Theorem \ref{ModRing2} (a). Clearly $P_S(k)$ is the number of partitions of $\frac{k}{q-1}$ with parts in $\{\frac{q-1}{q-1},\frac{q^2-1}{q-1},\ldots,\frac{q^r-1}{q-1}\}$, which by \cite[\TheoremCite 15.2]{Nathanson} has the asymptotic behaviour given in (b).
Assertion (c) is a direct consequence of Theorem \ref{ModRing3}.
Finally, by Theorem \ref{ModRing2} (c) the dimension in (d) is just the number of partitions of $k$ into $r$ summands, which is well-known to be $\binom{k-1}{r-1}$.
\end{Proof}

\begin{Rem}\label{DimensionFormulasRem}
\rm Taking invariants one may obtain similar dimension formulas for arbitrary arithmetic subgroups $\Gamma$ containing $\Gamma(t)$. In particular \cite[Thm.\;8.4]{Pink} gives an explicit formula when $\Gamma(t)<\Gamma<\Gamma_1(t)$. It seems an interesting problem to find a dimension formula in general.
\end{Rem}


\begin{Rem}\label{CohenMacaulay}
\rm Combining Theorem \ref{ModRing1} and \cite[Thm.\;1.7]{PinkSchieder} shows that $\CM_{*}(\Gamma(t))$ is a Cohen-Macaulay normal integral domain. By taking invariants, the argument in \cite[\S2]{PinkSchieder} shows the same for $\CM_{*}(\Gamma)$ whenever $\Gamma(t)<\Gamma<\Gamma_1(t)$. For $\Gamma=\GL_r(A)$ and $\SL_r(A)$ the same follows from the explicit description in Theorem \ref{ModRing2}. One may ask: Is this only a rare event for small level, or is it a general phenomenon?
\end{Rem}



\begin{center}
\rule{8cm}{0.01cm}
\end{center}

\begin{minipage}[t]{5cm}{\small
Department of Mathematical Sciences \\
University of Stellenbosch \\
Stellenbosch, 7600 \\
South Africa \\
djbasson@sun.ac.za
}
\end{minipage}\hfill
\begin{minipage}[t]{5cm}{\small
School of Mathematical and Physical Sciences \\
University of Newcastle \\
Callaghan, 2308 \\
Australia \\
florian.breuer@newcastle.edu.au\\
 {\em and} \\
Department of Mathematical Sciences \\
University of Stellenbosch \\
Stellenbosch, 7600 \\
South Africa \\}
\end{minipage}\hfill
%
%
\begin{minipage}[t]{5cm}{\small
Department of Mathematics \\
ETH Z\"urich\\
8092 Z\"urich\\
Switzerland \\
pink@math.ethz.ch}
\end{minipage}

\end{document}